\documentclass{amsart}
\usepackage{latexsym}
\usepackage{amssymb}
\usepackage{amsmath}
\usepackage{tikz}
\def\ni{\noindent}
\def\be{\begin{equation}}
\def\ee{\end{equation}}
\def\bC{\mathbb{C}}
\def\bR{\mathbb{R}}

\def\bN{\mathbb{N}}
\def\bZ{\mathbb{Z}}
\def\gg{\mathfrak{g}}
\def\gk{\mathfrak{k}}
\def\gp{\mathfrak{p}}
\def\gt{\mathfrak{t}}
\def\gsu{\mathfrak{su}}
\def\gso{\mathfrak{so}}
\def\gsp{\mathfrak{sp}}
\def\gspin{\mathfrak{spin}}

\def\B{\mathrm{B}}
\def\e{\mathrm{e}}
\def\h{\widehat}
\def\ol{\overline}
\def\dsp{\displaystyle}
\newtheorem{df}{Definition}[section]
\newtheorem{prop}[df]{Proposition}

\newtheorem{remark}[df]{Remark}

\newtheorem{lem}[df]{Lemma}

\begin{document}

\title[The First Eigenvalue of the Dirac Operator.../the outer case]
{The First Eigenvalue of the Dirac Operator on Compact Outer Spin Symmetric Spaces}%
\author{Jean-Louis Milhorat}
\address{Laboratoire Jean Leray\\
 UMR CNRS 6629\\
 D\'epartement de Math\'ematiques\\
 Universit\'e de Nantes\\
 2, rue de la Houssini\`ere\\
 BP 92208\\
  F-44322 NANTES CEDEX 03}
\email{jean-louis.milhorat@univ-nantes.fr}


\begin{abstract}
In two previous papers, we started a study of the first eigenvalue of the Dirac operator on compact spin symmetric spaces, 
providing, for symmetric spaces of ``{inner}'' type, a formula giving this first eigenvalue in terms of the algebraic data of the groups involved.
We conclude here that study by giving the explicit expression of the first eigenvalue for ``{outer}'' compact spin symmetric spaces.
\end{abstract}

\maketitle
\section{Introduction}
It is well-known that symmetric spaces provide examples where
the spectrum of Laplace or Dirac operators
can be (theoretically) explicitly computed. However this explicit computation is far from being simple
in general and only a few examples are known.  On the other hand, several classical results in geometry 
involve the first (nonzero) eigenvalue of those spectra, so it seems interesting to get this eigenvalue without
computing all the spectrum. In two previous papers (see \cite{Mil05} and \cite{Mil06}), we stated a formula 
giving the square of the first eigenvalue of the Dirac operator of a spin compact symmetric space in terms of the algebraic data 
of the groups involved. However, this formula was based on a result of R. Parthasarathy, \cite{Par}, only valid 
for symmetric spaces of inner type. Recall that a symmetric space $G/K$ is said to be of inner type if 
the involution characterizing it is given by an inner conjugation in the group $G$, or alternatively, by the fact 
that the groups $G$ and $K$ have same rank, i.e. own a common maximal torus. Otherwise, it is said to be outer,
(see for instance Sec.~8.6 in \cite{Wol} for details).
The study of subgroups $K$ of maximal rank in a compact Lie group was initiated by A. Borel and J. De Siebenthal
in \cite{BdS}, with an explicit description for compact simple groups. In \cite{Mur52} (see also \cite{Mur65}), S. Murakami gave a 
general method to study outer involutive automorphisms of compact simple Lie algebras. Using those results, 
the following complete list of 
irreducible compact simply-connected Riemannian symmetric spaces $G/K$ of type I with 
$\mathrm{rank}\, K<\mathrm{rank}\, G$, can be obtained (see \cite{Mur65} or
Sec.~8.12 in J. A. Wolf's book \cite{Wol})~:
\begin{gather*}
\mathrm{SU}(2m)/\mathrm{SO}(2m)\;;\; \mathrm{SU}(2m+1)/\mathrm{SO}(2m+1)\;;\; 
\mathrm{SU}(2m)/\mathrm{Sp}(m)\;;\\
\mathrm{SO}(2p+2q+2)/\mathrm{SO}(2p+1)\times\mathrm{SO}(2q+1)\;;\; 
\mathrm{E}_6/\mathrm{F}_4\;;\; \mathrm{E}_6/\mathrm{Sp}(4)\,.
\end{gather*}
It was proven by M. Cahen and S. Gutt in \cite{CG}, that all symmetric spaces in that list, except
$\mathrm{SU}(2m+1)/\mathrm{SO}(2m+1)$, are spin.

In the present paper, the following explicit value for the square of the first eigenvalue of the Dirac operator is 
obtained for all those symmetric spaces $G/K$, endowed with the Riemannian metric induced by the Killing
form of $G$ sign-changed.

\vspace{5mm}
\begin{center}
\begin{tabular}{|c|l|}
\hline
Symmetric space&\hfill Square of the first eigenvalue\hfill\hfill\\
&\hfill of the Dirac operator\hfill\hfill\\
\hline
&\\
& $m$ even,\; $\frac{1}{12}\,(m+1)(4m-1)+\frac{1}{32}$\\
$\displaystyle{\frac{\mathrm{SU}(2m)}{\mathrm{SO}(2m)}}$ &\\
&$m$ odd,\; $\frac{1}{12}\,(m+1)(4m-1)+\frac{1}{32}\left(1-\frac{1}{m^2}\right)$ \\
&\\
\hline &\\
&$m$ even,\; $\frac{1}{12}\,(m-1)(4m+1)+\frac{1}{32}$\\
$\displaystyle{\frac{\mathrm{SU}(2m)}{\mathrm{Sp}(m)}}$ &\\
&$m$ odd $\geq 3$,\; $\frac{1}{12}\,(m-1)(4m+1)+\frac{1}{32}\left(1-\frac{1}{m^2}\right)$ \\
&\\
\hline
&\\
$\frac{\mathrm{SO}(2p+2q+2)}{\mathrm{SO}(2p+1)\times\mathrm{SO}(2q+1)}$&$\frac{1}{16\,(p+q)}\, \big(
8pq\,(2p+q+1)+4p\,(p+1)+4q\,(q+1)+1\big)$\\
 \small{$p\leq q$}&\\
 &\\
\hline
&\\
$\displaystyle{\frac{\mathrm{E}_6}{\mathrm{F}_4}}$ & \hfill$\displaystyle{\frac{277}{72}}$\hfill\hfill\\
&\\
\hline
&\\
$\displaystyle{\frac{\mathrm{E}_6}{\mathrm{Sp}_4}}$ &\hfill $\displaystyle{\frac{529}{72}}$\hfill\hfill\\
&\\
\hline

\end{tabular}
\end{center}

\section{Preliminaries for the proof}

\subsection{Spectrum of the Dirac operator on spin compact irreducible symmetric spaces}
From now on, we consider a spin compact simply connected irreducible symmetric space $G/K$ of ``type I'',
where $G$ is a simple compact and simply-connected Lie group and $K$ is the
connected subgroup formed by the fixed elements of an involution
$\sigma$ of $G$.  This involution induces the
Cartan decomposition of the Lie algebra $\gg$ of $G$ into
$$\gg=\gk\oplus\gp\,,$$ where $\gk$ is the Lie algebra of $K$ and
$\gp$ is the vector space $\{X\in\gg\,;\, \sigma_{*}\cdot X=-X\}$.
The symmetric space $G/K$ is endowed with the Riemannian metric induced by the restriction to $\gp$ of
the Killing form $\B_G$ of $G$ sign-changed.

\noindent The spin condition implies that the homomorphism
$$\alpha: h\in K \longmapsto \mathrm{Ad}_{G}(h)_{|\gp}\in \mathrm{SO}(\gp)$$
lifts to a homomorphism $\widetilde{\alpha}:K\longrightarrow\mathrm{Spin}(\gp)$
such that $\xi\circ \widetilde{\alpha}=\alpha$ where $\xi$ is the two-fold covering
$\mathrm{Spin}(\gp)\rightarrow \mathrm{SO}(\gp)$, \cite{CG}.

\noindent Then the group $K$ inherits a spin representation given by
$$\widetilde{\rho_K}: K \xrightarrow{\widetilde{\alpha}} \mathrm{Spin}(\gp)\xrightarrow{\rho}
\mathrm{GL}_{\bC}{(\Sigma)}\,,$$
where $\rho$ is the spinor representation in the complex spinor space $\Sigma$.

\noindent The Dirac operator has a real discrete spectrum, symmetric with respect to the origin.
A real number $\lambda$ belongs to the spectrum if and only if there exists an irreducible representation
$\gamma:G\rightarrow \mathrm{GL}_{\bC}{(V_{\gamma})}$ whose restriction $\mathrm{Res}_K^G(\gamma)$ to the subgroup $K$,
contains in its decomposition into irreducible parts, a representation equivalent to some irreducible component
of the decomposition of the spin representation $\widetilde{\rho_K}$ of $K$.
Then
\begin{equation}
\label{eigenv}\lambda^2=c_{\gamma}+n/16\,,
\end{equation}
where $c_\gamma$ is the Casimir eigenvalue of the irreducible representation $\gamma$
(which only depends on the equivalence class of $\gamma$)
and where $n=\dim (G/K)$, $n/16$ being $\mathrm{Scal}/8$ for the choice of the metric (cf. \cite{BH3M} or
 \cite{Gin09} for details).
 
\ni Hence the first eigenvalue of the Dirac operator is given by the lowest $c_\gamma$, among the irreducible representations
$\gamma's$ of $G$ such that $\mathrm{Res}_K^G(\gamma)$ contains an irreducible component of the spin representation $\widetilde{\rho_K}$ of $K$.
We will say for short that such an irreducible representation $\gamma$ (or its highest weight) verifies the
``{spin condition}''.

\subsection{Outer symmetric spaces}

\ni Let $T$ be a maximal torus of $G$.  
Then $T_K:=T\cap K$ is a maximal torus of $K$. 

\ni As it was already mentioned, the symmetric space $G/K$ is said to be outer if the involution $\sigma: G\rightarrow G$ is not
a conjugation in the group, and this is equivalent to the condition  $\dim(T_K)<\dim(T)$.

\ni Let $\gt$ and $\gt_K$ be the Lie algebras of $T$ and $T_K$, and let $\gt_0:=\gt\cap \gp$. 
Note that $\gt_0$ is the orthogonal complement of $\gt_K$ in $\gt$ for 
the scalar product $-\B_G$.

\ni Let $\gg_{\bC}$ and $\gk_{\bC}$ be the complexifications of the Lie algebras $\gg$ and $\gk$, and $\langle\,\,,\,\rangle$ the $\bC$-extension of
$-\B_G$ to $\gg_{\bC}$. The root decompositions of $\gg_{\bC}$ and $\gk_{\bC}$ under the respective actions of $T$ and $T_K$ are given by
\begin{equation}\label{ggdecomp}
\gg_{\bC}=\gt_{\bC}\oplus\left(\mathop{\oplus}_{i=1}^{N} \gg_{\pm \theta_i}\right)\,,
\end{equation}
where $\pm\theta_i$ are the $G$-roots, and 
$$ \gk_{\bC}=\gt_{K,\bC}\oplus\left(\mathop{\oplus}_{i=1}^{M} \gk_{\pm \theta'_i}\right)\,,$$
where $\pm\theta'_i$ are the $K$-roots.

\ni The set $\Phi$ of $G$-roots decomposes into the disjoint union  
$$\Phi=\Phi_1\cup \Phi_2\cup \Phi_3\,,$$
where 
\begin{align*}
 \Phi_1&=\{\theta \in \Phi, \gg_{\theta}\subset \gk_{\bC}\}\,,\\
 \Phi_2&=\{\theta \in \Phi, \gg_{\theta}\subset \gp_{\bC}\}\,,\\
 \Phi_3&=\Phi\backslash (\Phi_1\cup \Phi_2)\,.
\end{align*}

\ni Note that if $\theta\in \Phi_i$, $i=1,2,3$, then $-\theta\in \Phi_i$, since $\gg_{-\theta}=\overline{\gg_{\theta}}$.

\ni We use the same notation for the involution $\sigma_{*}:\gg\rightarrow \gg$ and its $\bC$-linear extension to
$\gg_{\bC}$. As $\sigma_{*}(\gt)=\gt$, by means of the 
scalar product $-\B_{G}$, $\sigma_{*|\gt}$ induces an involution $\sigma^{*} : i\,\gt^{*}\rightarrow i\,\gt^{*}$. The scalar product 
on $i\,\gt^{*}$ induced by $-\B_{G}$  is denoted by $\langle\;_,\, \rangle$.

\ni We consider for any $G$-root $\theta_i$, a basis $E_{\theta_i}$ of the one-dimensional space $\gg_{\theta_i}$, and set
$$ U_{\theta_i}:=\frac{1}{2}\,\left(E_{\theta_i}+\sigma_{*}E_{\theta_i}\right)\;\text{and}\;
V_{\theta_i}:=\frac{1}{2}\,\left(E_{\theta_i}-\sigma_{*}E_{\theta_i}\right)\,.$$
Then $U_{\theta_i}\in \gk$, $V_{\theta_i}\in \gp$, $E_{\theta_i}=U_{\theta_i}+V_{\theta_i}$, and
$$\theta_i\in \Phi_3\Longleftrightarrow (U_{\theta_i}\not=0\; \text{and}\; V_{\theta_i}\not=0)\,.$$

\ni In order to determine an expression of the weights of the spin representation $\widetilde{\rho_K}: K\rightarrow\mathrm{GL}_{\bC}{(\Sigma)}$, 
we will have to use the following properties of $G$-roots. For any root $\theta\in \Phi$, the restriction $\theta_{|\gt_K}$ is 
denoted $\theta'$ for short.

\begin{lem}\label{lem.roots}
 \begin{description}
  \item[R1] $\forall \alpha\in \Phi_3$, $\exists X\in \gt_K$ such that $\alpha(X)\not =0$, hence $\alpha'$ is 
  a (non zero) $K$-root.
  \item[R2] $\forall \alpha\in \Phi$, $\alpha\in \Phi_1\cup \Phi_2\Longleftrightarrow \alpha_{|\gt_0}=0$.
  Hence if $\alpha$ and $\beta$ are two roots in $\Phi_1$ such that $\beta\not = \pm \alpha$, then
  $\alpha'$ and $\beta'$ are two $K$-roots such that $\beta'\not =\pm \alpha'$.
  \item[R3] If $\alpha\in \Phi_1$ and $\beta\in \Phi_3$, then $\alpha'$ and $\beta'$ are two $K$-roots 
  such that $\beta'\not=\pm \alpha'$.
  \item[R4] For all $\alpha\in \Phi$, $\sigma^{*}(\alpha)$ is a root such that 
  \begin{enumerate}
   \item if $\alpha\in \Phi_1\cup \Phi_2$, then $\sigma^{*}(\alpha)=\alpha$,
   \item if $\alpha\in \Phi_3$, then $\sigma^{*}(\alpha)\not=\pm\alpha$,
  \end{enumerate}
 \item[R5] For all $\alpha$ and $\beta\in \Phi_3$, $\beta'=\alpha'\Longleftrightarrow \beta=\alpha$ or 
 $\beta=\sigma^{*}(\alpha)$.
\end{description}
\end{lem}
\proof  The proofs are given in Appendix. Some of those results appear in some way in Chapter~3 of \cite{BR90}.\qed

\ni According to \textbf{R4}, we consider $\Phi_3'=\{\gamma_1,\ldots,\gamma_p\}\subset\Phi_3$ such that $\Phi_3$ is 
the disjoint union
$$\Phi_3=\{\gamma_1,\ldots,\gamma_p\}\cup \{\sigma^{*}(\gamma_1),\ldots, \sigma^{*}(\gamma_p)\}\,.$$
Then
\begin{lem} The set $\Phi_K$ of $K$-roots under the action of $T_K$ is given by
$$\Phi_K=\{\theta_{|\gt_K};\, \theta \in \Phi_1\cup \Phi_3'\}\,.$$
\end{lem}
\proof By definition, the restriction to $\gt_K$ of a root in $\Phi_1$ is a $K$-root, and by
\textbf{R1}, the restriction to $\gt_K$ of a root in $\Phi_3'$ is also a $K$-root.
Conversely, let $E'_{\varrho}\in \gk_{\bC}$ be a root-vector for a $K$-root $\varrho$.
As for any $\alpha\in \Phi_3$, $\sigma_{*}(E_\alpha)$ is a root-vector for the root $\sigma^{*}(\alpha)$
by \eqref{R11} and \eqref{R12}, we may write, according to \eqref{ggdecomp}, 
\begin{align*}E'_{\varrho}=H+\sum_{\alpha\in \Phi_1}\lambda_{\alpha}\, E_\alpha
+\sum_{\alpha\in \Phi_2}\lambda_{\alpha}\, E_\alpha &+\sum_{\alpha\in \Phi_3'}\lambda_{\alpha}\, E_\alpha
+\sum_{\alpha\in \Phi_3'}\lambda_{\sigma^{*}(\alpha)}\, \sigma_{*}(E_\alpha)\,,\\
& \qquad\text{where $H\in \gt$ and $\lambda_{\alpha}\in \bC$,}\\
=H+\sum_{\alpha\in \Phi_1}\lambda_{\alpha}\, E_\alpha+\sum_{\alpha\in \Phi_2}\lambda_{\alpha}\, E_\alpha 
&+ \sum_{\alpha\in \Phi_3'}\left(\lambda_{\alpha}+\lambda_{\sigma^{*}(\alpha)}\right)\, U_{\alpha}\\
&+ \sum_{\alpha\in \Phi_3'}\left(\lambda_{\alpha}-\lambda_{\sigma^{*}(\alpha)}\right)\, V_{\alpha}\,.
\end{align*}
Since $\sigma_{*}(E'_{\varrho})=E'_{\varrho}$, one gets $H\in \gt_K$, $\forall \alpha\in \Phi_2$, $\lambda_{\alpha}=0$, and
$\forall \alpha\in \Phi_3'$, $\lambda_{\sigma^{*}(\alpha)}=\lambda_{\alpha}$, hence
$$E'_{\varrho}=H+\sum_{\alpha\in \Phi_1}\lambda_{\alpha}\, E_\alpha+2\,\sum_{\alpha\in \Phi_3'}\lambda_{\alpha}\, U_\alpha\,.$$
Now, as for any $X\in \gt_K$, $[X,E'_{\varrho}]=\varrho(X)\, E'_{\varrho}$, one gets
\begin{align*}\forall X\in \gt_K\,,\; -\varrho(X)\, H&+\sum_{\alpha\in \Phi_1}\lambda_{\alpha} \left(\alpha(X)
 -\rho(X)\right)\, E_\alpha\\
&+2\,\sum_{\alpha\in \Phi_3'}\lambda_{\alpha} \left(\alpha(X)-\rho(X)\right)\, U_\alpha=0\,.\end{align*}
Hence $H=0$, and, as there exists at least one $\lambda_{\alpha}\not=0$, $\varrho=\alpha'$. Note furthermore
that by \textbf{R2}, \textbf{R3} and \textbf{R5}, there exists only one such a $\lambda_{\alpha}\not=0$. \qed

\ni By the results \textbf{R1} and \textbf{R2}, the restriction to $\gt_K$ of any $G$-root is nonzero.
As $\gt_K$ can not be the finite union of the hyperplanes $\mathrm{Ker}(\alpha')$ for $\alpha\in \Phi$, 
there exists $X\in \gt_K$ such that for any $\alpha\in \Phi$, $\alpha(X)\not=0$. This implies that $X$
is regular in $\gt$, and also regular in $\gt_K$, by the description of $K$-roots given in the above lemma. 
We define the set $\Phi^{+}$ (resp. $\Phi_K^{+}$) of positive roots of $G$ (resp. $K$) by the condition
$$\theta \in \Phi^{+}\;\text{(resp. $\Phi_K^{+}$)}\;\Longleftrightarrow \theta(X)>0\,.$$
Note that as $X\in \gt_K$, one has $\theta\in \Phi^{+}\Longleftrightarrow \sigma^{*}(\theta)\in \Phi^{+}$, 
and that by the above considerations
\begin{equation}\label{pos.K.roots}\Phi_K^{+}=\{\theta_{|\gt_K};\, \theta \in \Phi_1^{+}\cup \Phi_3'^{+}\}\,.
\end{equation}
Note furthermore that, by the description of positive $G$-roots, the half-sum $\delta_G$ of the positive roots verifies
$$\sigma^{*}(\delta_G)=\delta_G\,.$$
\subsection{Weights of the spin representation of $K$.}\label{w.spin}
The decomposition of the spin representation depends on the parity of $\dim \gp_{\bC}$, so we begin by considering
the even dimensional case.
\subsubsection{First case~: $\dim \gp_{\bC}$ even.}
As it will be seen below, this amounts to suppose that $\dim \gt_0$ is even. Setting $\dim \gt_0=2r_0$,
we consider an orthonormal basis $(T_k)$, $k=1,\ldots 2r_0$, of $\gt_0$,
and construct the Witt basis $(Z_k,\overline{Z_k})$, $k=1,\ldots r_0$, defined by 
$$Z_k:=\frac{1}{2}\, (T_{2k-1} +i\, T_{2k})\,,\quad \text{and}\quad \overline{Z_k}:=\frac{1}{2}\, (T_{2k-1} -i\, T_{2k})\,.$$
For the $\bC$-linear extension $\langle\;,\,\rangle$ of the scalar product on $\gg$, one has
$$\langle Z_i, Z_j\rangle =\langle \overline{Z_i}, \overline{Z_j}\rangle=0\,,\quad \text{and}\quad 
\langle Z_i, \overline{Z_j}\rangle=\frac{1}{2}\, \delta_{ij}\,.$$
Let $\Phi_1^{+}=\{\alpha_1,\ldots,\alpha_\ell\}$, $\Phi_2^{+}=\{\beta_1,\ldots,\beta_p\}$ and 
$\Phi_3'^{+}=\{\gamma_1,\ldots,\gamma_q \}$. 

\begin{lem}\label{root-vect.basis} Vectors $E_\theta$, $E_{-\theta}$, $\theta=\beta_1,\ldots,\beta_p$, 
$V_\theta$, $V_{-\theta}$, $\theta=\gamma_1\,\ldots,\gamma_q$, may be choosen such that, with 
$Z_i$, $\overline{Z_i}$, $i=1,\ldots, r_0$, they define a Witt basis of $\gp_{\bC}$, in the sense that
\begin{align*}
\langle E_\alpha,E_\beta\rangle =\langle E_{-\alpha},E_{-\beta}\rangle=
\langle E_\alpha,V_\beta\rangle =\langle E_{-\alpha},V_{-\beta}\rangle
&=\langle V_\alpha,V_\beta\rangle =\langle V_{-\alpha},V_{-\beta}\rangle\,,\\
&=\langle Z_i, Z_j\rangle =\langle \overline{Z_i}, \overline{Z_j}\rangle=0\,,
\end{align*}
and
\begin{equation}\label{Wbasis}
\langle E_\alpha, E_{-\beta}\rangle=\frac{1}{2}\,\delta_{\alpha \beta}\,,\quad
\langle V_\alpha, V_{-\beta}\rangle=\frac{1}{2}\,\delta_{\alpha \beta}\,,\quad
\langle Z_i, \overline{Z_j}\rangle=\frac{1}{2}\, \delta_{ij}\,.
\end{equation}
\end{lem}
\proof Clearly, any such vectors span $\gp_{\bC}$: any $X\in \gp_{\bC}$ may be written as  
\begin{align*}
 H&+\sum_{\alpha\in \Phi_1}\lambda_{\alpha}\, E_\alpha+\sum_{\alpha\in \Phi_2}\lambda_{\alpha}\, E_\alpha
+ \sum_{\alpha\in \Phi_3'}\left(\lambda_{\alpha}+\lambda_{\sigma^{*}(\alpha)}\right)\, U_{\alpha}\\
&+ \sum_{\alpha\in \Phi_3'}\left(\lambda_{\alpha}-\lambda_{\sigma^{*}(\alpha)}\right)\, V_{\alpha}\,,
\end{align*}
with $H\in \gt_{\bC}$ and $\lambda_\alpha\in \bC$.
The condition $\sigma_{*}(X)=-X$ implies $H\in \gt_{0,\bC}$ and $\forall \alpha\in \Phi_1$, $\lambda_\alpha=0$,
and $\forall \alpha\in \Phi_3'$, $\lambda_\alpha=-\lambda_{\sigma^{*}(\alpha)}$, hence
$$X=H+\sum_{\alpha\in \Phi_2^{+}} \lambda_\alpha\, E_\alpha+\sum_{\alpha\in \Phi_2^{+}} \lambda_{-\alpha}\,
E_{-\alpha}+2\, \sum_{\alpha\in \Phi_3'^{+}} \lambda_\alpha\, V_\alpha
+2\, \sum_{\alpha\in \Phi_3'^{+}} \lambda_{-\alpha}\, V_{-\alpha}\,,$$
so $X$ is a linear combination of the considered vectors.

\ni On the other hand, it is well-known that, for any couple of roots $(\alpha,\beta)$ such that $\alpha+\beta\not=0$,
$\langle \gg_\alpha,\gg_\beta\rangle=0$. Hence
\begin{itemize}
 \item if $\alpha$ and $\beta\in \Phi_2$ are such that $\alpha+\beta\not=0$, then
 $\langle E_\alpha,E_\beta\rangle =\langle E_{-\alpha},E_{-\beta}\rangle=0$.
 \item if $\alpha\in \Phi_2$ and $\beta\in \Phi_3$, then from $\langle E_\alpha,E_\beta\rangle =0$, one obtains as
 $\langle  E_\alpha, U_\beta\rangle =0$, $ \langle E_\alpha, V_\beta\rangle =0$, and then
 $\langle  E_{-\alpha}, V_{-\beta}\rangle =0$.
 \item if $\alpha\in \Phi_3'$ and $\beta\in \Phi_3'$ are such that $\alpha+\beta\not=0$, then $\alpha'+\beta'\not =0$ (otherwise by \textbf{R5}, 
 $\beta=\sigma^{*}(-\alpha)$, which is impossible), so $\langle  U_\alpha, U_\beta\rangle =0$, and then
 $\langle  V_\alpha, V_\beta\rangle =0$, (so $\langle  V_{-\alpha},V_{-\beta}\rangle =0$ also).
\end{itemize}
Finally, as $\langle\;,\,\rangle$ is non-degenerate, $(E_\alpha,E_{-\alpha})\not = 0$, $\alpha\in \Phi_2^{+}$, 
and $(V_\alpha,V_{-\alpha})\not = 0$, $\alpha\in \Phi_3'^{+}$, and we may choose the basis $(E_\alpha)$ such that 
\eqref{Wbasis} is verified. All of those orthogonality relations imply that the considered vectors are linearly independent.
\qed

\ni Since $E_{-\alpha}$ belongs to $\ol{\gg_{\alpha}}$, we may suppose that 
$E_{-\beta_{j}}=\ol{E_{\beta_{j}}}$ and $V_{-\gamma_{j}}=\ol{V_{\gamma_{j}}}$.

\ni In the Clifford algebra $\bC \ell(\gp)$, let 
$$\ol{w}:=\ol{E_{\beta_1}}\cdots \ol{E_{\beta_p}}\cdot
\ol{V_{\gamma_1}}\cdots \ol{V_{\gamma_q}}\cdot \overline{Z_{1}}\cdots \overline{Z_{r_0}}\,.$$
For any $I=\{i_1,\ldots,i_a\}\subset \{1,\ldots,p\}$, $i_1<i_2<\cdots <i_a$, 
(resp. $J=\{j_1,\ldots,j_b\}\subset \{1,\ldots,q\}$, $j_1<j_2<\cdots <j_b$;
resp. $K=\{k_1,\ldots,k_c\}\subset \{1,\ldots,r_0\}$, $k_1<k_2<\cdots <k_c$), we introduce the notation
$$E_I\cdot V_J\cdot Z_K:=E_{\beta_{i_1}}\cdots E_{\beta_{i_a} }\cdot V_{\gamma_{j_1}}\cdots V_{\gamma_{j_b} }\cdot
Z_{k_1}\cdots Z_{k_c}\,,$$
setting $E_I$ (resp. $V_J$, resp. $Z_K$) $=1$, if $I$ (resp. $J$, resp. $K$) $=\emptyset$.
As $I$ (resp. $J$, resp. $K$) runs through the set of subsets of $\{1,\ldots,p\}$, (resp.
$\{1,\ldots,q\}$, resp. $\{1,\ldots,r_0\}$), the vectors $E_I\cdot V_J\cdot Z_K\cdot \ol{w}$ define a basis of the spinor
space $\Sigma:=\bC \ell(\gp)\cdot \ol{w}$ (cf. \cite{BH3M}).

\ni Considering
\begin{equation}\label{onbp}\begin{cases} X_j&=E_{\beta_{j}}+\ol{E_{\beta_{j}}}\,,\\
Y_j&=i\,\left(E_{\beta_{j}}-\ol{E_{\beta_{j}}}\right)\,,\\
& 1\leq j\leq p\,,\end{cases}\quad \text{and}\quad
\begin{cases} F_j&=V_{\gamma_{j}}+\ol{V_{\gamma_{j}}}\,,\\
G_j&=i\,\left(V_{\gamma_{j}}-\ol{V_{\gamma_{j}}}\right)\,,\\
& 1\leq j\leq q\,,\end{cases}\end{equation}
one obtains an orthonormal basis $(X_i,Y_i,F_j,Y_j,T_k)_{\substack{1\leq i\leq p,\\
1\leq j\leq q,\\ 1\leq k\leq 2r_0}}$ of $\gp$.

\ni Then for any $X\in \gt_K$, one  has 
$$\begin{cases}[X,X_j]&=-i\, \beta_j(X)\, Y_j\,,\\
[X,Y_j]&=i\, \beta_j(X)\, X_j\,,\\
& 1\leq j \leq p\,,\end{cases}\; ;\;
\begin{cases}[X,F_j]&=-i\, \gamma_j(X)\, G_j\,,\\
[X,G_j]&=i\, \gamma_j(X)\, F_j\,,\\
& 1\leq j \leq q\,,\end{cases}\;\text{and}\;
\begin{cases}
 [X,T_j]=0\,,\\
 1\leq j\leq 2r_0\,.
\end{cases}$$
So
$$\forall X\in \gt_K\,,\quad\alpha_{*}(X)=-i\,\sum_{j=1}^{p}\beta_j(X)\, X_j\wedge Y_j
-i\,\sum_{j=1}^{q}\gamma_j(X)\, F_j\wedge G_j\,,$$
hence
$$\forall X\in \gt_K\,,\quad\widetilde{\alpha}_{*}(X)=-\frac{1}{2}\,i\,\sum_{j=1}^{p}\beta_j(X)\, X_j\cdot Y_j
-\frac{1}{2}\,i\,\sum_{j=1}^{q}\gamma_j(X)\, F_j\cdot G_j\,.$$

\ni Now it is easy to verify that
\begin{align*}
 X_k\cdot Y_k\cdot E_I\cdot V_J\cdot Z_K\cdot \ol{w}&=\begin{cases} -i\, E_I\cdot V_J\cdot Z_K\cdot \ol{w}\quad
 \text{if $k\notin I$,}\\
 i\, E_I\cdot V_J\cdot Z_K\cdot \ol{w}\quad
 \text{if $k\in I$,}\end{cases}\\
 F_k\cdot G_k\cdot E_I\cdot V_J\cdot Z_K\cdot \ol{w}&=\begin{cases} -i\, E_I\cdot V_J\cdot Z_K\cdot \ol{w}\quad
 \text{if $k\notin J$,}\\
 i\, E_I\cdot V_J\cdot Z_K\cdot \ol{w}\quad
 \text{if $k\in J$.}\end{cases}
 \end{align*}
 Hence 
 \begin{multline}
 \forall X\in \gt_K\,,\quad 
 \widetilde{\rho_{K}}_{*}(X)(E_I\cdot V_J\cdot Z_K\cdot \ol{w})=\nonumber\\
 \frac{1}{2}\left(\sum_{i\in I} \beta_{i}(X)-\sum_{i\notin I} \beta_{i}(X)+
 \sum_{j\in J} \gamma_{j}(X)-\sum_{j\notin J} \gamma_{j}(X)\right)\, E_I\cdot V_J\cdot Z_K\cdot \ol{w}\,.
 \end{multline}
So the $E_I\cdot V_J\cdot Z_K\cdot \ol{w}$ are weight-vectors, and a $\mu \in i\, \gt_K^{*}$ is a weight 
if and only if it can be expressed as 
\begin{equation}\label{weight}
 \frac{1}{2}\left( \pm \beta_1'\pm\cdots\pm\beta_p'\pm \gamma_1'\pm\cdots\pm\gamma_q'\right)\,.
\end{equation}
Note that as $K$ runs through the set of subsets of $\{1,\ldots,r_0\}$, all of the 
$E_I\cdot V_J\cdot Z_K\cdot \ol{w}$ for a given $I$ and $J$, belong to same weight-space,
hence the multiplicity of a weight of the form \eqref{weight} is $2^{r_0}$ times the number of ways in which
it can be expressed in the given form.

\ni Considering the volume element 
$$\omega:=i^{p+q+r_0}\, X_1\cdot Y_1\cdots X_p\cdot Y_p\cdot F_1\cdot G_1\cdots F_q\cdot G_q\cdot
T_1\cdot T_2\cdots T_{2 r_0-1}\cdot T_{2 r_0}\,,$$
one has 
$$\omega\cdot E_I\cdot V_J\cdot Z_K\cdot \ol{w}=(-1)^{\# I+\# J+\# K}\,E_I\cdot V_J\cdot Z_K\cdot \ol{w}\,,
$$
hence the spinor space $\Sigma$ decomposes into two irreducible components 
$$\Sigma=\Sigma^{+}\oplus \Sigma^{-}\,,$$
where $\Sigma^{+}$ (resp. $\Sigma^{-}$) is the linear span of the $E_I\cdot V_J\cdot Z_K\cdot \ol{w}$ 's 
such that $\# I+\# J+\# K$ is even (resp. odd).

\subsubsection{Second case~: $\dim \gp_{\bC}$ odd.}\label{odd.case}
This corresponds to the case $\dim \gt_0$ odd. Setting $n=\dim \gp$, by the choice of an orthonormal basis such as \eqref{onbp},
$\mathrm{SO}(\gp)$ is identified with $\mathrm{SO}(n)$, which is itself embedded in $\mathrm{SO}(n+1)$ in such a 
way that $\mathrm{SO}(n)$ acts trivially on the last vector $e_{n+1}$ of the standard basis of $\bR^{n+1}$. 

\ni Setting $\dim \gt_0=2r_0-1$, $r_0\geq 1$, the spinor space 
can be described (cf. \cite{BH3M}), as the space of positive spinors introduced above in the even-dimensional case :
$\Sigma=\mathrm{Span}\{E_I\cdot V_J\cdot Z_K\cdot \ol{w},\;\# I+\# J+\# K\;\text{even}\}$. By the result above,
the weights may be expressed as \eqref{weight}, but now, the multiplicity of such a weight is
$2^{r_0-1}$ times the number of ways in which
it can be expressed in the given form.

\subsubsection{A remark on highest weights}
\begin{lem} \label{h.weight} Any highest weight of the spinor representation of $K$ has necessarily the form
\begin{equation}
 \frac{1}{2}\left( \pm \beta_1'\pm\cdots\pm\beta_p'+ \gamma_1'+\cdots+\gamma_q'\right)\,.
\end{equation}
\end{lem}
\proof Let $\mu$ be a weight such that for a $j_0\in \{1,\ldots,q\}$, $\gamma_{j_0}'$ appears with 
a minus sign in the expression of $\mu$. Then $\mu+\gamma_{j_0}'$ is also a weight, hence $\mu$ can not be dominant,
since $\gamma_{j_0}'$ is a positive $K$-root. \qed

\section{The symmetric space $\dsp{\frac{\mathrm{SU}(2m)}{\mathrm{SO}(2m)}}$.}
In the following, $M_{2m}(\bC)$ denotes the space of 
$2m\times 2m$ matrices with complex coefficients, and $(E_{ij})$ its standard basis.
Here 
$$G=\mathrm{SU}(2m)=\{A\in M_{2m}(\bC), {}^t\overline{A}A=I_{2m}\}\,,$$
  and $$K=\{A\in \mathrm{SU}(2m)\;;\; \overline{A}=A\}=\mathrm{SO}(2m)\,.$$

\ni The symmetric space structure is given by the involution $\sigma:G\rightarrow G$,
$A\mapsto {}^t A^{-1}$. This involution induces the decomposition of the Lie algebra $\gsu_{2m}$
 of $\mathrm{SU}(2m)$ into:
$$\gsu_{2m}=\gso_{2m}\oplus \gp\,,$$
where $\gso_{2m}$ is the Lie algebra $\{X\in \gsu_{2m}\,;\, {}^t X=-X\}$ of $\mathrm{SO}(2m)$, and 
$$\gp=\{X\in \gsu_{2m}\,;\, {}^t X=X\}\,.$$
Elements of $\gp$ are traceless symmetric $2m\times 2m$ matrices with coefficients in $i\, \bR$, hence 
$\dim (\gp)= 2m^2+m-1$.

\ni Let $T_{s}=\left\{ \begin{pmatrix}\e^{i\theta_1}&&\\&\ddots&\\&& \e^{i\theta_{2m}}\end{pmatrix}
\;;\; \theta_i\in \bR\;,\; \sum_{i=1}^{2m} \theta_i=0\right\}$ be the standard maximal torus of
$\mathrm{SU}(2m)$. We consider a conjugate $T$ of $T_{s}$ in $G$ in such a way that $T\cap K$ is the 
standard maximal torus of $K$. For that, let
 $A_0\in \mathrm{SU}(2m)$ be the $2\times 2$ block diagonal matrix defined by
 $$A_0=\frac{\sqrt{2}}{2}\, \begin{pmatrix}\Box&&\\
                           &\ddots&\\
                            &&\Box
                          \end{pmatrix}\,, \quad\text{where}\; \Box=\begin{pmatrix} 1&i\\i&1\end{pmatrix}\,.
 $$
Then for any $H=\begin{pmatrix}\e^{i\theta_1}&&\\&\ddots&\\&& \e^{i\theta_{2m}}\end{pmatrix}\in T_{s}$,
$A_0 H A_0^{-1}=\begin{pmatrix} \Box_1&&\\
&\ddots&\\ &&\Box_m\end{pmatrix}$, where
$$\Box_j=
\e^{\frac{1}{2}i\,(\theta_{2j-1}+\theta_{2j})}\,\begin{pmatrix}\cos\left(\frac{1}{2}(\theta_{2j-1}-\theta_{2j})\right)&
\sin\left(\frac{1}{2}(\theta_{2j-1}-\theta_{2j})\right)\\-\sin\left(\frac{1}{2}(\theta_{2j-1}-\theta_{2j})\right)&
\cos\left(\frac{1}{2}(\theta_{2j-1}-\theta_{2j})\right)\end{pmatrix}\,.$$
Hence 
$A_0 H A_0^{-1}\in K \Longleftrightarrow A_0 H A_0^{-1}=\begin{pmatrix} \Box_1'&&\\
&\ddots&\\ &&\Box_m'\end{pmatrix}$, where
$$
\Box_j'=\begin{pmatrix}\cos(\theta_{2j})&-\sin(\theta_{2j})\\
\sin(\theta_{2j})&\cos(\theta_{2j})\end{pmatrix}\,,$$
so $A_0 H A_0^{-1}\in K$ if and only if it belongs to the standard maximal torus $T_K$ of $K$.

\ni The Lie algebra $\gt$ of $T$ is defined by
\begin{align*}\gt=\Big\{\sum_{j=1}^{m} i\gamma_{2j-1}\, (E_{2j-1\, 2j-1}+E_{2j\, 2j})&+\gamma_{2j}\,
(E_{2j\, 2j-1}-E_{2j-1\, 2j})\;;\;\\
&\qquad\qquad \gamma_{j}\in \bR\,,\sum_{j=1}^{m} \gamma_{2j-1}=0
\Big\}\,,\end{align*}
and 
\begin{align*}
\gt_K&=\left\{\sum_{j=1}^{m} \gamma_{2j}\,(E_{2j\, 2j-1}-E_{2j-1\, 2j})\,;\,\gamma_{2j}\in
\bR\right\}\,,\\
\gt_0&=\left\{\sum_{j=1}^{m} i\gamma_{2j-1}\, (E_{2j-1\, 2j-1}+E_{2j\, 2j})\,,\,\gamma_{2j-1}\in
\bR\,,\,\sum_{j=1}^{m} \gamma_{2j-1}=0\right\}\,. 
\end{align*}

\ni It is useful to consider the elements  $\h{x}_{j}$, $j=1,\ldots, 2m$ of $i\,\gt^{*}$,  defined by  considering the restriction to $\gt$ of 
the dual basis of the family of vectors $i\,(E_{2j-1\, 2j-1}+E_{2j\, 2j})$,  $(E_{2j\, 2j-1}-E_{2j-1\, 2j})$, 
$j=1,\ldots,m$~:
\begin{equation}\begin{split}
 \forall H&=\sum_{j=1}^{m} i\gamma_{2j-1}\, (E_{2j-1\, 2j-1}+E_{2j\, 2j})+\gamma_{2j}\,
(E_{2j\, 2j-1}-E_{2j-1\, 2j})\in\gt\,,\nonumber\\
&\h{x}_{2j-1}(H)=i\, \gamma_{2j-1}\quad \text{and}\quad \h{x}_{2j}(H)=i\,\gamma_{2j}\,.
\end{split}\end{equation}
The scalar product on $i\,\gt^{*}$ is given by the scalar product on $\gsu_{2m}$ defined by
$$\langle X,Y\rangle =-\frac{1}{2}\, \Re\big(\mathrm{Tr}(XY)\big)=-\frac{1}{8m} \mathrm{B}(X,Y)\,,\quad
X, Y\in \gsu_{2m}\,.$$

\begin{lem}
 Each $\mu\in i\,\gt^{*}$ may be uniquely written as 
 \begin{equation}\label{w}
 \mu=\sum_{i=1}^{2m} \mu_i\, \h{x}_i\,,\; \text{where $\mu_i\in\bR$ verify $\sum_{i=1}^{m} \mu_{2i-1}=0$,}
 \end{equation}
 and, for any $\mu=\sum_{i=1}^{2m} \mu_i\, \h{x}_i$ and $\mu'=\sum_{i=1}^{2m} \mu_i'\, \h{x}_i \in i\, \gt^{*}$,
 $$\langle \mu,\mu'\rangle=\sum_{i=1}^{2m} \mu_i\,\mu_i'\,.
 $$
\end{lem}
\proof As the $\h{x}_{2j-1}$, $j=1,\ldots, m-1$ and $\h{x}_{2j}$, $1,\ldots, m$ define a basis of $\gt^{*}$, any 
$\mu\in\gt^{*}$ is uniquely written as $\mu=\sum_{i=1}^{m-1} \mu_{2i-1}\, \h{x}_{2i-1}+\sum_{i=1}^{m}\mu_{2i}\,
\h{x}_{2j}$, $\mu_i\in\bR$. Now since $\h{x}_{2m-1}=-\sum_{i=1}^{m-1} \h{x}_{2i-1}$, setting 
$\mu_{2m-1}=-\frac{1}{m}\,\sum_{i=1}^{m-1}\mu_{2j-1}$, one gets
$$\mu =\sum_{i=1}^{m-1} (\mu_{2i-1}+\mu_{2m-1})\, \h{x}_{2i-1}+\mu_{2m-1}\, \h{x}_{2m-1}+\sum_{i=1}^{m}\mu_{2i}\,
\h{x}_{2i}\,,$$
so $\mu$ may be written as \eqref{w}. Now if $\mu=\sum_{i=1}^{2m} \lambda_i\, \h{x}_i$, where 
$\lambda_i\in\bR$ verify $\sum_{i=1}^{m} \lambda_{2i-1}=0$, then 
$\mu=\sum_{i=1}^{m-1} (\lambda_{2i-1}-\lambda_{2m-1})\, \h{x}_{2i-1}+
\sum_{i=1}^{m} \lambda_{2i}\, \h{x}_{2i}$, and so 
$$\mu_{2i-1}=\lambda_{2i-1}-\lambda_{2m-1}\,,\;i=1,\ldots, m-1\,,\quad\text{and}\quad
\mu_{2i}=\lambda_{2i}\,,\;i=1,\ldots, m\,.$$ Then $\mu_{2m-1}=-\frac{1}{m}\,\sum_{i=1}^{m-1}\mu_{2i-1}=\lambda_{2m-1}$,
hence $\mu_{2i-1}+\mu_{2m-1}=\mu_{2i-1}+\lambda_{2m-1}=\lambda_{2i-1}$, hence the unicity of the writting \eqref{w}.

\ni The last result follows from the fact that the vectors $E_{2j-1\, 2j-1}+E_{2j\, 2j}$,
 $E_{2j\, 2j-1}-E_{2j-1\, 2j}$, $j=1,\ldots, m$, define an orthonormal set for the scalar product $\langle\;,\,\rangle$. \qed 

\ni In the following, any $\mu\in i\, \gt^{*}$, of the form \eqref{w} will be denoted
$$\mu=[\mu_1,\mu_3,\ldots,\mu_{2m-1};\mu_2,\mu_4,\ldots,\mu_{2m}]\,,\quad \sum_{i=1}^{m} \mu_{2i-1}=0\,.$$
Note that the $\h{x}_{2i}'$, $1\leq i\leq m$, define a basis of $i\,\gt_K^{*}$. Any $\lambda\in i\,\gt_K^{*}$ of the form $\lambda=\sum_{i=1}^{m} \lambda_{2i}\, \h{x}_{2i}'$, 
$\lambda_{2i}\in \bR$, will be denoted 
$$\lambda=(\lambda_2,\lambda_4,\ldots, \lambda_{2m})\,.$$
The involution $\sigma^{*}$ of $i\gt^{*}$ induced by $\sigma$ is defined by
$$[\mu_1,\mu_3,\ldots,\mu_{2m-1};\mu_2,\mu_4,\ldots,\mu_{2m}]\xrightarrow{\sigma^{*}} [-\mu_1,-\mu_3,\ldots,-\mu_{2m-1};\mu_2,\mu_4,\ldots,\mu_{2m}]\,.$$

\subsection{Sets of roots.}
Since the root-vectors relative to the standard torus $T_{s}$ are the $E_{ij}$, the root-vectors relative 
to the torus $T$ are given by the $A_0 E_{ij} A_0^{-1}$. Explicitly, the roots are
\begin{align*}
&(\h{x}_{2i-1}-\h{x}_{2j-1})+(\h{x}_{2i}-\h{x}_{2j})\,,\; 1\leq i\not =j\leq m\,, \quad\text{root-space~:
$\bC\,A_0 E_{2i-1\, 2j-1} A_0^{-1}$,}\\
 &(\h{x}_{2i-1}-\h{x}_{2j-1})-(\h{x}_{2i}-\h{x}_{2j})\,,\; 1\leq i\not =j\leq m\,, \quad\text{root-space~: 
 $\bC\, A_0 E_{2i\, 2j} A_0^{-1}$,}\\
&(\h{x}_{2i-1}-\h{x}_{2j-1})+(\h{x}_{2i}+\h{x}_{2j})\,,\; 1\leq i,j\leq m\,, \quad\text{root-space~:
$\bC\,A_0 E_{2i-1\, 2j} A_0^{-1}$,}\\
&(\h{x}_{2i-1}-\h{x}_{2j-1})-(\h{x}_{2i}+\h{x}_{2j})\,,\; 1\leq i,j\leq m\,, \quad\text{root-space~:
$\bC\,A_0 E_{2i\, 2j-1} A_0^{-1}$.}
 \end{align*}
Note that the only roots that verify $\sigma^{*}(\theta)=\theta$ are $\pm 2\, \h{x}_{2i}$, and, as they are not
$K$-roots, one may conclude from lemma~\ref{lem.roots} that
$$ \Phi_1=\emptyset\,,\quad \Phi_{2}=\{\pm 2\, \h{x}_{2i}\,,\, i=1,\ldots,m\}\,,$$
and all other roots belong to $\Phi_3$.

\ni We choose positive roots such that
\begin{align*}
 \Phi_2^{+}&=\{2\, \h{x}_{2i}\;,\; i=1,\ldots,m\}\,,\\
 \Phi_3'^{+}&=\left\{\begin{array}{cc}
 (\h{x}_{2i-1}-\h{x}_{2j-1})+(\h{x}_{2i}-\h{x}_{2j})\,,& 1\leq i<j\leq m\,,\\
 (\h{x}_{2i-1}-\h{x}_{2j-1})+(\h{x}_{2i}+\h{x}_{2j})\,,& 1\leq i<j\leq m\,,\\
 \end{array}\right\}\,,\\
 \sigma^{*}(\Phi_3'^{+})&=\left\{\begin{array}{cc}
 -(\h{x}_{2i-1}-\h{x}_{2j-1})+(\h{x}_{2i}-\h{x}_{2j})\,,& 1\leq i<j\leq m\,,\\
 -(\h{x}_{2i-1}-\h{x}_{2j-1})+(\h{x}_{2i}+\h{x}_{2j})\,,& 1\leq i<j\leq m\,\\
 \end{array}\right\}\,.
\end{align*}
It is easy to verify that a system of simple roots is given by
\begin{align*}
 \beta:=&2\, \h{x}_{2m}\,,\\
 \gamma_{i}:= &(\h{x}_{2i-1}-\h{x}_{2i+1})+(\h{x}_{2i}-\h{x}_{2(i+1)})\,,\quad i=1,\ldots, m-1\,,\\
 \sigma^{*}(\gamma_{i})=&-(\h{x}_{2i-1}-\h{x}_{2i+1})+(\h{x}_{2i}-\h{x}_{2(i+1)})\,,\quad i=1,\ldots, m-1\,.
 \end{align*}
Note that 
\begin{align*}
 \langle \gamma_i, \gamma_i\rangle&=\langle \sigma^{*}(\gamma_i), \sigma^{*}(\gamma_i)\rangle=
 \langle \beta,\beta\rangle=4\,,\quad i=1,\ldots, m-1\,,\\
 \langle \gamma_i, \gamma_{i+1}\rangle&=\langle \sigma^{*}(\gamma_i), \sigma^{*}(\gamma_{i+1})\rangle=-2\,,\; 
 i=1,\ldots,m-2\,,\\
 \langle \gamma_{m-1},\beta\rangle&=\langle \sigma^{*}(\gamma_{m-1}),\beta\rangle= -2,\, \,,
 \end{align*}
all the other scalar products being zero. Hence
\begin{lem} \label{dom.w.G} A vector $\mu=[\mu_1,\mu_3,\ldots,\mu_{2m-1};\mu_2,\mu_4,\ldots,\mu_{2m}]\in i\,\gt^{*}$,
$\sum_{i=1}^{m} \mu_{2i-1}=0$, is a dominant weight if and only if 
\begin{align*}
 &\left.\begin{array}{cc}&\mu_{2i}-\mu_{2(i+1)}\in \bN\,,\\
 &\mu_{2i-1}-\mu_{2i+1}\in \bZ\,,\end{array}\right\}\;\text{both simultaneously odd or even,}\\
 & |\mu_{2i-1}-\mu_{2i+1}|\leq \mu_{2i}-\mu_{2(i+1)}\,, \quad 1\leq i\leq m-1\,,\\
 &\mu_{2m}\in \bN\,.
\end{align*}
\end{lem}
\proof Note just that the conditions 
$2\, \frac{\langle \gamma_{i}, \mu\rangle}{\langle \gamma_{i},\gamma_{i}\rangle}\in
\bN$, $2\, \frac{\langle \sigma^{*}(\gamma_{i}), \mu\rangle}{\langle \gamma_{i},\gamma_{i}\rangle}\in
\bN$, $i=1,\ldots m-1$, are equivalent to
\begin{align*}&\frac{(\mu_{2i-1}-\mu_{2i+1})+(\mu_{2i}-\mu_{2(i+1)})}{2}\in \bN\,,\\
\intertext{and}\;
&\frac{-(\mu_{2i-1}-\mu_{2i+1})+(\mu_{2i}-\mu_{2(i+1)})}{2}\in \bN\,,
\end{align*}
which imply (and are equivalent to) the first three conditions of the lemma. \qed

\ni As it was remarked before, the positive $K$-roots are given by considering the restrictions 
to $\gt_K$ of the positive roots in $\Phi_3'$, so 
$$\Phi_K^{+}=\left\{ \h{x}_{2i}'-\h{x}_{2j}'\,,\, \h{x}_{2i}'+\h{x}_{2j}'\,,\,
1\leq i<j\leq m\right\}\,.$$
A set of simple roots is given by
\begin{align*}
 \theta_{i}'&= \h{x}_{2i}'-\h{x}_{2(i+1)}'\,,\quad i=1,\ldots, m-1\,,\\
 \theta_{m}'&= \h{x}_{2(m-1)}'-\h{x}_{2m}'\,.
\end{align*}
Note that 
\begin{align*}
 \langle \theta_{i}',\theta_{i}'\rangle &=2\,,\quad \langle \theta_{i}',\theta_{i+1}'\rangle =-1\,,\quad
 1\leq i\leq m-2\,,\\
 \langle \theta_{m-2}',\theta_{m}'\rangle &=-1\,,\quad \langle \theta_{m-1}',\theta_{m}'\rangle=0\,,
\end{align*}
hence the following ``{classical}'' characterization~:
\begin{lem} \label{dom.w.K} A vector $\lambda=(\lambda_2,\lambda_4,\ldots, \lambda_{2m})\in i\,\gt_K^{*}$,
is a dominant weight if and only if 
 $$\lambda_2\geq \lambda_4\geq \cdots \geq \lambda_{2(m-1)}\geq |\lambda_{2m}|\,,$$
 and the $\lambda_{2i}$ are all simultaneously integers or half-integers.
\end{lem}
\subsection{Highest weights of the spin representation of $K$.}
\begin{prop}
 The spin representation of $K$ has two highest weights:
$$(m,m-1,\ldots,2,\pm 1)\,,$$
 both with multiplicity $2^{[\frac{m-1}{2}]}$.
\end{prop}
\proof One has $\frac{1}{2}\, \sum_{\gamma\in \Phi_3'^{+} } \gamma'=\sum_{1\leq i<j\leq m} \h{x}_{2i}'=
\sum_{i=1}^{m-1} (m-i)\,\h{x}_{2j}'$. Hence by the result of lemma~\ref{h.weight}, any highest weight $\lambda$ has 
necessarily the form
$$\lambda= \sum_{i=1}^{m} (m-i\pm 1)\, \h{x}_{2i}'\,.$$
But the dominance condition of lemma~\ref{dom.w.K} implies then 
$$\lambda=\sum_{i=1}^{m-1} (m-i+1)\, \h{x}_{2i}'\pm\, \h{x}_{2m}'=(m,m-1,\ldots,2,\pm 1)\,.$$
Denote by $\lambda_{\pm}$ those two dominant weights.
Let $\delta_{K}$ be the half-sum of the positive $K$-roots. One has 
$\delta_K=(m-1,m-2,\ldots,1,0)$, hence
$$\langle \lambda_{+}+\delta_K, \lambda_{+}+\delta_K \rangle = \langle \lambda_{-}+\delta_K, \lambda_{-}+\delta_K \rangle\,.$$
Hence the two weights are both highest weights, since otherwise, one of the two should be contained in the set
of weights of an irreducible representation having the other one as highest weight, 
and the above equality should be impossible (cf. lemma~C, 13.4 in \cite{Hum}).

\ni With the help of the Weyl dimension formula, it may be checked that any irreducible module with 
highest weight $\lambda_{\pm}$ has dimension $2^{(m-1)(m+1)}$. As we noticed it before, the multiplicity of each one 
of those two weights is at least $2^{\frac{m-1}{2}}$ if $m$ is odd, and $2^{\frac{m-2}{2}}$ if $m$ is even.
Since if $m$ is odd,
$$ 2\times 2^{\frac{m-1}{2}}\times 2^{(m-1)(m+1)}=2^{\frac{2m^2+m-1}{2}}=\dim(\Sigma)\,,$$
and if $m$ is even,
$$2\times 2^{\frac{m-2}{2}}\times 2^{(m-1)(m+1)}=2^{\frac{2m^2+m-2}{2}}=\dim(\Sigma)\,,$$
one concludes that the multiplicity of each weight is exactly $2^{\left[\frac{m-1}{2}\right]}$. \qed

\subsection{The first eigenvalue of the Dirac operator.}
Recall that the first eigenvalue of the Dirac operator is given by the lowest Casimir eigenvalue
$c_\gamma$, among the irreducible representations $\gamma$ of $G$ verifying the ``{spin condition}''~:
$\mathrm{Res}_K^G(\gamma)$ contains an irreducible component of the
spin representation $\widetilde{\rho_K}$ of $K$. By the Freudenthal formula, the Casimir eigenvalue of 
an irreducible representation with highest weight $\mu_{\gamma}$ is given by
$$c_\gamma= \langle \mu_{\gamma},\mu_{\gamma}+2\, \delta_G\rangle=\|\mu_{\gamma}+\delta_G\|^2 -\|\delta_G\|^2\,.$$
Hence we look for $G$-dominant weights $\mu_{\gamma}$, verifying the spin condition  and such that 
$\|\mu_{\gamma}+\delta_G\|^2$ is minimal.

\ni We first determine the $G$-weights $\mu$ (non necessarily dominant) for which
$\mu_{|\gt_K}=\lambda_{\pm}$ and $\|\mu+\delta_G\|^2$ is minimal.

\ni Note first that 
$\mu=[\mu_1,\mu_3,\ldots,\mu_{2m-1};\mu_{2},\mu_{4},\ldots,\mu_{2m}]$, $\sum_{i=1}^{m}
\mu_{2i-1}=0$, is a $G$-weight if and only if  $\mu_{2i-1}-\mu_{2i+1}\in \bZ$, $i=1,\ldots m-1$,
$\mu_{2i}\in \bZ$, $i=1,\ldots, m$, and  $\mu_{2i-1}-\mu_{2i+1}$, $\mu_{2i}-\mu_{2(i+1)}$
are both simultaneously even or odd, $i=1,\ldots m-1$.

\ni Let $\mu$ be such a $G$-weight. It verifies $\mu_{|\gt_K}=\lambda_{\pm}$ if and only if it has the form
$$\mu=[\mu_1,\mu_3,\ldots,\mu_{2m-1}; m,m-1,\ldots,2,\pm 1]\,,$$
where $\sum_{i=1}^{m}\mu_{2i-1}=0$, and $\mu_{2i-1}-\mu_{2i+1}$ are odd integers,
$i=1,\ldots m-1$.

\ni The half-sum of the positive $G$-roots $\delta_G$ is given by
\begin{align*}\delta_G&=2\sum_{1\leq k<l\leq m} \h{x}_{2k}+\sum_{k=1}^{m}\h{x}_{2k}=2\, \sum_{k=1}^{m}
\left(m-k+\frac{1}{2}\right)\, \h{x}_{2k}\,,\\
&= [0,0,\ldots, 0; (2m-1),(2m-3),\ldots, 1]\,.
\end{align*}
Hence $\|\mu+\delta_G\|^2$ is minimal if and only if $\sum_{i=1}^{m} \mu_{2i-1}^2$ is minimal.

\ni Set $\mu_{2i-1}-\mu_{2i+1}=k_i$, $i=1,\ldots, m-1$.
For any $j=1,\ldots,m-1$, one has
\begin{align*}p_j:=\sum_{i=1}^{j}k_i&=\mu_1+\sum_{i=1}^{j-1}\mu_{2i+1}-\sum_{i=1}^{j-1}\mu_{2i+1}-\mu_{2j+1}\,,\\
 &=\mu_1-\mu_{2j+1}\,.
\end{align*}
Using $\sum_{i=1}^{m}\mu_{2i-1}=0$, one then gets
$$\sum_{j=1}^{m-1}p_j=(m-1)\, \mu_1+\mu_1=m\, \mu_1\,,$$
hence 
$$\mu_1=\frac{1}{m}\sum_{i=1}^{m-1}p_i\,,$$
and
$$\mu_{2j+1}=\mu_1-p_j=\frac{1}{m}\left(\sum_{i=1}^{m-1} \,p_i-m\,p_j\right)\,.$$
The expression $\sum_{i=1}^{m} \mu_{2i-1}^2$ is a polynomial $F(p_1,\ldots,p_{m-1})$ of the variables 
$p_1$, $p_2\ldots,p_{m-1}$. With the notation $p=(p_1,\ldots,p_{m-1})$, one has
$$\frac{\partial F}{\partial p_i}(p)=-2\,\mu_{2i+1}\,.$$
With no surprise, $F$ as only one critical point at $(0,\ldots,0)$.

\ni Now $$\frac{\partial^2 F}{\partial p_i\partial p_j}(0)=-\frac{2}{m}\,(1-m\,\delta_{ij})\,.$$
Denote by $H$ the Hessian matrix of $F$ at $0$. It has for eigenvalues $\frac{2}{m}$ with multiplicity one 
and $2$ with multiplicity $m-2$, and the following vectors define an orthogonal basis of eigenvectors.
\begin{align*}
 &v_1= (1,1,\ldots,1)\,,\\
 &v_i= (0,0,\ldots,0,m-i,\underbrace{-1,-1,\ldots,-1}_{m-i})\,,\; 2\leq i\leq m-1\,.
\end{align*}
Considering the orthonormal basis $(\epsilon_i\, v_i)_{1\leq i\leq m-1}$, $\epsilon_i=1/\|v_i\|$, and 
denoting by $Q$ the orthogonal matrix of the change of basis, one gets
\begin{align*}
 F(p)&=\frac{1}{2}\, {}^t p H p\,,\\
 &=\frac{1}{m}\, {}^t p\, Q\begin{pmatrix} 1 &0& 0&0\\
 0&m&0&0\\ &&\ddots&\\0&0&0&m\end{pmatrix}{}^t Q\,p\,,\\
 &=\frac{1}{m}\,\epsilon_1\,\big((m-1)k_1+(m-2)k_2+\cdots +k_{m-1}\big)^2\\
 &\qquad + \epsilon_2\, \big((m-2)k_2+(m-3)k_3\cdots +k_{m-1}\big)^2\\
 &\qquad +\cdots\\
 &\qquad + \epsilon_{m-i}\, \big(i\,k_{m-i}+\cdots +k_{m-1}\big)^2\\
 &\qquad +\cdots\\
 &\qquad +\epsilon_{m-2}\, \big(2 k_{m-2}+k_{m-1}\big)^2\\
 &\qquad +\epsilon_{m-1}\, \big(k_{m-1}\big)^2\,.
\end{align*}
The minimum is obtained if and only if all the squares are minimal. Hence, since the $k_i$ have to be odd integers, 
\begin{lem}\label{min}
 The minimum is obtained only when 
 $$(k_1,k_2,\ldots,k_{m-1})=(1,-1,1,-1,\ldots)\;\text{or}\;
 (-1,1,-1,1,\ldots)\,.$$
\end{lem}
\proof  Assume that all the squares are minimal. First, since $k_{m-1}$ is an odd integer, 
$k_{m-1}^2$ is minimum only when $k_{m-1}=\pm 1$. Now, let us show by
induction that $k_{m-i}=(-1)^{i+1}\, k_{m-1}$, $i=1,\ldots,m-1$.
Assuming that this is true for $1\leq i\leq 2j-1$, $j\geq 1$, one gets
\begin{align*}
 \big(2j\,k_{m-2j}+(2j-1)\, k_{m-(2j-1)}+\cdots+k_{m-1}\big)^2
 &=\big( 2j\,k_{m-2j}+j\,k_{m-1}\big)^2\,,\\
 &= j^2\,\big(2\, k_{m-2j}+k_{m-1}\big)^2\,.
\end{align*}
Since $2\, k_{m-2j}+k_{m-1}$ is an odd integer, the above square is minimal if and only if
$\big(2\, k_{m-2j}+k_{m-1}\big)^2=1$, which implies
$$ 4\, k_{m-2j}\,(k_{m-2j}+k_{m-1})=0\,,$$
hence $k_{m-2j}=-k_{m+1}=(-1)^{2j+1}\, k_{m-1}$. Thus the result is also true for $i=2j$, and one gets
$$\big((2j+1)\,k_{m-(2j+1)}+2j\, k_{m-(2j-1)}+\cdots+k_{m-1}\big)^2
=\big( (2j+1)\,k_{m-(2j+1)}-j\,k_{m-1}\big)^2\,.$$
If $k_{m-1}=1$, then $(2j+1)\,k_{m-(2j+1)}-j\,k_{m-1}\geq j+1$ if $k_{m-(2j+1)}\geq 1$ and
$(2j+1)\,k_{m-(2j+1)}-j\,k_{m-1}\leq -(3j+1)$ if $k_{m-(2j+1)}\leq -1$, thus the above square is 
minimal if and only if $k_{m-(2j+1)}= 1=k_{m-1}=(-1)^{2j+2}\, k_{m-1}$.

\ni If $k_{m-1}=-1$, then $(2j+1)\,k_{m-(2j+1)}-j\,k_{m-1}\geq 3j+1$ if $k_{m-(2j+1)}\geq 1$ and
$(2j+1)\,k_{m-(2j+1)}-j\,k_{m-1}\leq -(j+1)$ if $k_{m-(2j+1)}\leq -1$, thus the above square is 
minimal if and only if $k_{m-(2j+1)}=-1=k_{m-1}=(-1)^{2j+2}\, k_{m-1}$. \qed

\ni Hence the $G$-weights $\mu$ for which $\mu_{|\gt_K}=\lambda_{\pm}$ and $\|\mu+\delta_G\|^2$ is minimal are given by
\begin{enumerate}
 \item If $m$ is even $m=2p$, 
 \begin{enumerate} 
 \item if $(k_1,k_2,\ldots,k_{m-1})=(1,-1,1,-1,\ldots,-1,1)$, then since $p_{2i-1}=1$, $i=1,\ldots,p$, and $p_{2i}=0$,
 $i=1,\ldots,p-1$,
 $$\mu_{\pm}=\left[\frac{1}{2},-\frac{1}{2},\ldots,\frac{1}{2},-\frac{1}{2}; m,m-1,\ldots,2,\pm 1\right]\,, $$
 \item if $(k_1,k_2,\ldots,k_{m-1})=(-1,1,-1,1,\ldots,1,-1)$, then since $p_{2i-1}=-1$, $i=1,\ldots,p$, and $p_{2i}=0$,
 $i=1,\ldots,p-1$,
 $$\mu_{\pm}'=\left[-\frac{1}{2},\frac{1}{2},\ldots,-\frac{1}{2},\frac{1}{2}; m,m-1,\ldots,2,\pm 1\right]\,. $$ 
 Note that $\mu_{\pm}'=\sigma^{*}(\mu_{\pm})$ and  $\|\mu_{\pm}+\delta_G\|^2=\|\mu_{\pm}'+\delta_G\|^2$.

 \end{enumerate}
 \item If $m$ is odd, $m=2p+1$,
 \begin{enumerate}
 \item if $(k_1,k_2,\ldots,k_{m-1})=(1,-1,1,-1,\ldots,1,-1)$, then since $p_{2i-1}=1$ and $p_{2i}=0$,
 $i=1,\ldots,p$,
$$\mu_{\pm}=\left[\frac{1}{2}-\frac{1}{2m},-\frac{1}{2}-\frac{1}{2m}\ldots,-\frac{1}{2}-\frac{1}{2m},\frac{1}{2}-\frac{1}{2m}; m,m-1,\ldots,2,\pm 1\right]\,,$$
\item if $(k_1,k_2,\ldots,k_{m-1})=(-1,1,-1,1,\ldots,-1,1)$, then since $p_{2i-1}=-1$ and $p_{2i}=0$,
 $i=1,\ldots,p$,
 $$\mu_{\pm}'=\left[-\frac{1}{2}+\frac{1}{2m},\frac{1}{2}+\frac{1}{2m},\ldots,\frac{1}{2}+\frac{1}{2m},-\frac{1}{2}+\frac{1}{2m}; m,m-1,\ldots,2,\pm 1\right]\,,$$
Here also note that $\mu_{\pm}'=\sigma^{*}(\mu_{\pm})$ and  $\|\mu_{\pm}+\delta_G\|^2=\|\mu_{\pm}'+\delta_G\|^2$. 
 \end{enumerate}
\end{enumerate}
Note that, by lemma~\ref{dom.w.G}, the weights $\mu_{+}$ and $\mu_{+}'$ are $G$-dominant, whereas 
$\mu_{-}$ (resp. $\mu_{-}'$) belongs to the orbit of $\mu_{+}$ (resp. $\mu_{+}'$) under the Weyl group, since 
$\mu_{-}=\sigma_\beta(\mu_{+})$, (resp. $\mu_{-}'=\sigma_\beta(\mu_{+}')$), where $\sigma_\beta$ is the reflexion
across $\beta^{\bot}$, $\beta$ being the simple root $2\,\h{x}_{2m}$.

\ni In order to conclude, we first remark that the $G$-dominant weights $\mu_{+}$ and 
$\mu_{-}$ verify the spin condition. That follows from the following general result.

\begin{lem} Let $\mu_{\gamma}$ be the highest weight of an irreducible representation $\gamma$ of $G$.
Then $\mu_{\gamma|\gt_K}$ is a $K$-dominant weight, and any irreducible representation of $K$ having
$\mu_{\gamma|\gt_K}$ as a highest weight is contained in the restriction $\mathrm{Res}_K^G(\gamma)$ of $\gamma$ to $K$.
\end{lem}
\proof By the result of lemma~\ref{dom.w.G}, a dominant $G$-weight
$$\mu=[\mu_1,\mu_3,\ldots,\mu_{2m-1};\mu_2,\mu_4,\ldots,\mu_{2m}]\,,\quad \sum_{i=1}^{m} \mu_{2i-1}=0\,,$$
restricts to $\gt_K$ as $\mu_{|\gt_K}=\sum_{i=1}^{m} \mu_{2i}\, \h{x}_{2i}=(\mu_{2},\mu_{4},\ldots,\mu_{2m})$,
where the $\mu_{2i}$ are non-negative integers verifying the condition $\mu_{2}\geq \mu_{4}\geq \cdots \geq \mu_{2m}$,
hence, by the result of lemma~\ref{dom.w.K}, $\mu_{|\gt_K}$ is a $K$-dominant weight.

\ni Now let $v_{\gamma}$ be the\footnote{it is unique up to a scalar multiple.} maximal vector of the representation $\gamma$. Since it is killed by the action
of root-vectors corresponding to positive roots, it is in particular killed by the action of the $E_{\alpha}$ and
$E_{\sigma^{*}(\alpha)}$, $\alpha\in {\Phi'_{3}}^{+}$, hence by the action of the $U_{\alpha}=\frac{1}{2}(E_{\alpha}
+\sigma_{*}(E_\alpha))$, $\alpha\in {\Phi'_{3}}^{+}$. Since the $U_\alpha$ are root-vectors for the positive
$K$-roots $\alpha'$, $\alpha\in {\Phi'_{3}}^{+}$, (see \eqref{pos.K.roots}), $v_\gamma$ is a maximal vector
of $\mathrm{Res}_K^G(\gamma)$ for the weight $\mu_{\gamma|\gt_K}$. \qed  

\ni The conclusion now results from the following remark.

\begin{lem} Let $\mu$ be a $G$-dominant weight. Let
$\Pi_\mu$ be the set of weights of any irreducible $G$-representation with highest weight $\mu$.  

\ni If $\lambda_{\pm}$ is the restriction to $\gt_K$ of a weight $\lambda\in \Pi_\gamma$ and 
if $\|\mu+\delta_G\|^2$ is minimal, then $\lambda$ belongs to the orbit of $\mu$ under 
the Weyl group $W_G$. 
\end{lem}
\proof The weight $\lambda$ lies in the orbit under the Weyl group $W_G$ of $G$ of a dominant weight
$\varrho\in \Pi_\gamma$. By the result of lemma~B, 13.3 in \cite{Hum}, $\|\lambda+\delta_G\|^2\leq \|\varrho+
\delta_G\|^2$, with equality only if $\lambda=\varrho$, and by lemma~C, 13.4 in \cite{Hum},
$\|\varrho+\delta_G\|^2\leq \|\mu+\delta_G\|^2$, with equality only if $\varrho=\mu$. Now as $\varrho$ is
a $G$-dominant weight which posses in its orbit a weight $\lambda$ whose restriction to $\gt_K$ is 
$\lambda_{\pm}$, the minimality condition verified by $\mu$ on that sort of weights implies 
$\|\varrho+\delta_G\|^2=\|\mu+\delta_G\|^2$, hence $\varrho=\mu$, and $\lambda$ belongs to the 
orbit of $\mu$. \qed

\begin{lem}\label{low.cas} If $\mu_\gamma$ is a $G$-dominant weight verifying the spin condition, then
$$\|\mu_\gamma+\delta_G\|^2\geq \|\mu_{+}+\delta_G\|^2\,.$$
\end{lem}
\proof As $\mu_\gamma$ verifies the spin condition, there exists a weight $\lambda\in \Pi_\gamma$ such 
that $\lambda_{|\gt_K}=\lambda_{\pm}$. By the result of the above lemma, $\lambda$ belongs to the orbit of
$\mu_\gamma$ under 
the Weyl group $W_G$, so $\|\lambda+\delta_G\|^2\leq \|\mu_\gamma+\delta_G\|^2$, with equality
if and only if $\lambda=\mu_\gamma$. If $\lambda=\mu_\gamma$, then, by the above considerations 
$\mu_\gamma =\mu_{+}$ or $\mu_{+}'$, and the result follows since $\|\mu_{+}+\delta_G\|^2=
\|\mu_{+}'+\delta_G\|^2$. If $\lambda\not =\mu_\gamma$, then by replacing $\lambda$ by 
$\sigma_\beta(\lambda)$ if necessary, we may suppose that $\lambda_{|\gt_K}=\lambda_{+}$. But then, by the above 
considerations, $\|\mu_{+}+\delta_G\|^2\leq \|\lambda+\delta_G\|^2$, and so 
$\|\mu_{+}+\delta_G\|^2<\|\mu_\gamma+\delta_G\|^2$.\qed

\ni Finally, we may conclude that the square of the first eigenvalue of the Dirac operator is given by 
$$
\frac{1}{8m}\, \langle \mu_{+},\mu_{+}+2\delta_G\rangle+\frac{2m^2+m-1}{16}\,,
$$
hence the result.

\section{The symmetric space $\dsp{\frac{\mathrm{SU}(2m)}{\mathrm{Sp}(m)}}$.}
Let $J$ be the matrix in $\mathrm{SU}(2m)$ defined by
$$J:=\begin{pmatrix} 0&-I_m\\I_m&0\end{pmatrix}\,.$$
Note that $J^2=-I_{2m}$ and $J^{-1}=-J={}^t J$.

\ni The group $\mathrm{Sp}(m)$ is identified with the subgroup of $\mathrm{SU}(2m)$ 
$$\mathrm{Sp}(m):=\{A\in \mathrm{SU}(2m)\,;\, {}^tA\,J\,A=J\}\,.$$

\ni The symmetric space structure is given by the involution 
$$\sigma:G\rightarrow G\,,\quad
A\mapsto {}^tJ\,{}^t A^{-1}\,J\,.$$ This involution induces the decomposition of the Lie algebra $\gsu_{2m}$
into:
$$\gsu_{2m}=\gsp_{m}\oplus \gp\,,$$
where $\gsp_{2m}$ is the Lie algebra 
$$\{X\in \gsu_{2m}\,;\, {}^t XJ=-JX\}=\{X\in \gsu_{2m}\,;\,
JXJ={}^tX\}\,,$$ of $\mathrm{Sp}(m)$, and 
$$\gp=\{X\in \gsu_{2m}\,;\, {}^t XJ=JX\}=\{X\in \gsu_{2m}\,;\,
JXJ=-{}^tX\}\,.$$
Elements of $\gp$ are matrices $X$ of the form $X=\begin{pmatrix} X_1&\ol{X_2}\\
X_2&-\ol{X_1}\end{pmatrix}$ where $X_1$ and $X_2$ are $m\times m$ matrices verifying
${}^t \ol{X_1}=-X_1$ and ${}^t X_2=-X_2$.
 hence 
$\dim (\gp)= 2m^2-m-1$.
Denoting by $T$ the maximal standard torus of $\mathrm{SU}(2m)$, $T\cap K$ is the standard maximal
torus of $\mathrm{Sp}(m)$~:
$$T_K:=T\cap K=\left\{\begin{pmatrix}\e^{i\beta_1}&&&&&\\&\ddots&&&&\\&&\e^{i\beta_m}&&&\\
&&&&&\\
&&&\e^{-i\beta_1}&&\\&&&&\ddots&\\&&&&&\e^{-i\beta_m}\end{pmatrix}\;;\; \beta_j\in \bR
\right\}\,.$$
The Lie algebra $\gt$ of $T$ is defined by
$$\gt=\left\{\sum_{j=1}^{2m} i \beta_j\, E_{jj}\;;\; \beta_j\in\bR\,,\;\sum_{j=1}^{2m}\beta_j=0\right\}\,,$$ 
and
\begin{align*}
 \gt_K&=\left\{\sum_{j=1}^{m} i \beta_j\, (E_{jj}-E_{m+j\,m+j})\;;\; \beta_j\in\bR\right\}\,,\\
 \gt_0&=\left\{\sum_{j=1}^{m} i \gamma_j\, (E_{jj}+E_{m+j\,m+j})\;;\; \gamma_j\in\bR\,,\;
 \sum_{j=1}^{m}\gamma_j=0\right\}\,.
\end{align*}
Let $\h{y}_j$, $1\leq j\leq 2m$, be the vectors of $i\, \gt^{*}$ defined by considering the restriction to $\gt$ of the dual basis of the family
of vectors $i\, E_{jj}$, $1\leq j\leq 2m$~:
$$\forall H=\sum_{j=1}^{2m} i\beta_j\, E_{jj}\in \gt\,,\quad \h{y}_j(H)=\beta_j\,.$$
Any element $\mu\in i\, \gt^{*}$ may be uniquely written as
$$\mu=\sum_{j=1}^{2m} \mu_{j}\,\h{y}_j\,,\quad \sum_{j=1}^{2m} \mu_{j}=0\,,$$
and is denoted 
$$\mu=[\mu_1,\ldots,\mu_{2m}]\,.$$
The scalar product on $i\,\gt^{*}$ considered here is given by the scalar product on $\gsu_{2m}$ defined by
$$\langle X,Y\rangle =- \Re\big(\mathrm{Tr}(XY)\big)=-\frac{1}{4m} \mathrm{B}(X,Y)\,,\quad
X, Y\in \gsu_{2m}\,.$$
For $\mu=[\mu_1,\ldots,\mu_{2m}]$ and $\mu'=[\mu_1',\ldots,\mu_{2m}']\in i\,\gt^{*}$,
$$\langle \mu,\mu'\rangle=\sum_{j=1}^{2m} \mu_{j}\,\mu_{j}'\,.$$
The involution $\sigma^{*}$ of $i\gt^{*}$ induced by $\sigma$ is defined by
$$\mu=[\mu_1,\ldots,\mu_{2m}]\xrightarrow{\sigma^{*}} [-\mu_{m+1},-\mu_{m+2},\ldots,-\mu_{2m},-\mu_{1},-\mu_{2},\ldots,-\mu_{m}]\,.$$
Note that the vectors $\h{x}_j=\h{y}_j':=\h{y}_{j|\gt_K}$, $j=1,\ldots,m$, define a basis of $i\gt_K^{*}$. Any $\lambda \in i\gt_K^{*}$ 
of the form $\lambda=\sum_{j=1}^{m}\lambda_j\, \h{x}_j$ is denoted 
$$\lambda=(\lambda_1,\ldots,\lambda_m)\,.$$

\subsection{Sets of roots.}
The root-vectors of $\mathrm{SU}(2m)$ relative to the standard torus are the $E_{ij}$. Explicitly the roots are 
\begin{align*}
 \h{y}_k-\h{y}_l\,,\quad & 1\leq k\not =l\leq m\,,\quad \text{with root-vector $E_{kl}$,}\\
 \h{y}_{m+k}-\h{y}_{m+l}\,,\quad & 1\leq k\not =l\leq m\,,\quad \text{with root-vector $E_{m+k\, m+l}$,}\\
 \pm (\h{y}_{k}-\h{y}_{m+l})\,,\quad & 1\leq k,l\leq m\,,\quad \text{with root-vector $E_{k\, m+l}$ (resp. $E_{m+l\, k}$).}
 \end{align*}
Note that the only roots $\theta$ that verify $\sigma^{*}(\theta)=\theta$ are $\pm (\h{y}_{k}-\h{y}_{m+k})$, and as 
$E_{k\, m+k}$ and $E_{m+k\, k}$ belong to $\gsp_{m}\otimes \bC$, those roots belong to $\Phi_1$. Hence
$$\Phi_1=\left\{\pm(\h{y}_{k}-\h{y}_{m+k})\,,\; 1\leq k\leq m\right\}\,,\quad \Phi_2=\emptyset\,,$$
and all other roots belong to $\Phi_3$.

\ni We choose positive roots such that
\begin{align*}
 \Phi_1^{+}&=\{\h{y}_{k}-\h{y}_{m+k}\,,\; 1\leq k\leq m\}\,,\\
 \Phi_{3}'^{+}&=\left\{\begin{array}{l} \h{y}_k-\h{y}_l\\ \h{y}_k-\h{y}_{m+l}\end{array}\,,\; 1\leq k<l\leq m\right\}\,,\\
 \sigma^{*}\big(\Phi_{3}'^{+}\big)&= \left\{\begin{array}{l} \sigma^{*}(\h{y}_k-\h{y}_l)=-\h{y}_{m+k}+\h{y}_{m+l}\\
 \sigma^{*}(\h{y}_k-\h{y}_{m+l})=-\h{y}_{m+k}+\h{y}_{l}\end{array}\,,\; 1\leq k<l\leq m\right\}\,.
\end{align*}
It is easy to verify that a system of simple roots is given by
\begin{align*}
 \alpha&:= \h{y}_m-\h{y}_{2m}\,,\\
 \gamma_k&:= \h{y}_k -\h{y}_{k+1}\,,\quad 1\leq k \leq m-1\,,\\
 \sigma^{*}(\gamma_k)&=-\h{y}_{m+k}+\h{y}_{m+k+1}\,,\quad 1\leq k \leq m-1\,.
\end{align*}
Note that 
\begin{align*}
 \langle \alpha,\alpha\rangle&=\langle \gamma_k, \gamma_k\rangle=\langle \sigma^{*}(\gamma_k),
 \sigma^{*}(\gamma_k)\rangle=2\,,\quad k=1,\ldots, m-1\,,\\
 \langle \gamma_k, \gamma_{k+1}\rangle&=\langle \sigma^{*}(\gamma_k), \sigma^{*}(\gamma_{k+1})\rangle=-1\,,
 \; k=1,\ldots,m-2\,,\\
 \langle \gamma_{m-1},\alpha\rangle&=\langle \sigma^{*}(\gamma_{m-1}),\alpha\rangle= -1\,,
 \end{align*}
all the other scalar products being zero. Hence
\begin{lem} \label{dom.w.G.2} 
A vector $\mu=[\mu_1,\ldots,\mu_{2m}]\in i\,\gt^{*}$,
$\sum_{j=1}^{2m} \mu_{j}=0$, is a dominant weight if and only if 
\begin{align*}
 &\left.\begin{array}{l}\mu_{k}-\mu_{k+1}\in \bN\,,\\
 \mu_{m+k+1}-\mu_{m+k}\in \bN\,,\end{array}\right\}\;k=1,\ldots,m-1\,,\\
 &\quad\mu_{m}-\mu_{2m}\in \bN\,.
\end{align*}
\end{lem}
\ni By \eqref{pos.K.roots}, the positive $K$-roots are the restrictions to $\gt_K$ of the positive $G$-roots in
$\Phi_1$ and $\Phi_3'$, hence, as $\h{y}_{m+k|\gt_K}=-\h{y}_{k|\gt_K}$,
$$\Phi_K^{+}=\left\{\begin{array}{ll} 2\, \h{x}_k \,, & 1\leq k\leq m\\ \h{x}_k -\h{x}_l\,, & 1\leq k<l\leq m\\
\h{x}_k +\h{x}_l\,, & 1\leq k<l\leq m\end{array}
\right\}\,.$$
A system of simple roots is given by
$$ \theta'_k:= \h{x}_k-\h{x}_{k+1}\,,\; k=1,\ldots, m-1\,,\quad \text{and}\quad \theta'_m:=2\,\h{x}_m\,.$$
Note that, as $\h{x}_k=\frac{1}{2}(\h{y}_k-\h{y}_{m+k})$, 
\begin{align*}
 \langle \theta'_k,\theta'_k\rangle &= 1\,,\quad k=1,\ldots, m-1\,, \quad \langle \theta'_m,\theta'_m\rangle = 2\,,\\
 \langle \theta'_k,\theta'_{k+1}\rangle &= -1/2\,,\quad k=1,\ldots, m-2\,,
 \quad \langle \theta'_{m-1},\theta'_m\rangle = -1\,,
\end{align*}
all the other scalar products being zero. Hence
\begin{lem} \label{dom.w.K.2} A vector $\lambda=(\lambda_1,\ldots \lambda_m)\in i\,\gt_K^{*}$,
is a dominant weight if and only if 
 $$\lambda_1\geq \lambda_2\geq \cdots \geq \lambda_{m}\geq 0\,,$$
 and all the $\lambda_{i}$ are integers.
\end{lem}

\subsection{Highest weights of the spin representation of $K$.}  

\ni Since $\Phi_2=\emptyset$, we may conclude from lemma~\ref{h.weight}~:
\begin{lem} The spin representation of $K$ has only one highest weight 
$$\lambda_0:=(m-1,m-2,\ldots,1,0)\,,$$
with multiplicity $2^{[\frac{m-1}{2}]}$
\end{lem}
\proof As we saw it before the multiplicity of $\lambda_0$ is at least $2^{[\frac{m-1}{2}]}$. But it may be checked,
using the Weyl dimension formula, that any irreducible $K$-module
with highest weight $\lambda_0$ has dimension $2^{m(m-1)}$, and
$$2^{[\frac{m-1}{2}]}\times 2^{m(m-1)}=\dim(\Sigma)\,,$$
hence the result.\qed

\subsection{The first eigenvalue of the Dirac operator.}

\ni We apply the same method as above, and begin by determining the $G$-weights $\mu$ such that
$\mu_{|\gt_K}=\lambda_0$ and $\|\mu+\delta_G\|^2$ is minimal.

\ni First a vector $\mu=[\mu_1,\ldots,\mu_{2m}]\in i\gt^{*}$, $\sum_{k=1}^{2m}\mu_k=0$, is a $G$-weight
if and only if $\mu_k-\mu_{k+1}\in \bZ$, $k=1,\ldots, 2m-1$.

\ni Such a vector $\mu$ verifies $\mu_{|\gt_K}=\lambda_0$ if and only if 
$$\mu_k-\mu_{m+k}=m-k\,,\quad k=1,\ldots, m\,.$$
The condition $\sum_{k=1}^{2m}\mu_k=0$ then implies that
$$\sum_{k=1}^{m}\mu_k=\frac{m(m-1)}{4}\,.$$
The half-sum of the positive $G$-roots is given by
$$\delta_G=\sum_{k=1}^m \left(m-k+\frac{1}{2}\right)\, (\h{y}_k -\h{y}_{m+k})\,.$$
(Note that $\sigma^{*}(\delta_G)=\delta_G$). Hence
$$\mu+\delta_G=\sum_{k=1}^{m}\left(\mu_k+m-k+\frac{1}{2}\right)\, \h{y}_k
+\sum_{k=1}^{m}\left(\mu_k-2(m-k)-\frac{1}{2}\right)\, \h{y}_{m+k}\,.$$
As $\sum_{k=1}^{m}\mu_k=\frac{m(m-1)}{4}$, $\|\mu+\delta_G\|^2$ is minimal if and only if
$\sum_{k=1}^{m} (\mu_k^2 +k\, \mu_k)$ is minimal.

\ni Set $\mu_{i}-\mu_{i+1}=k_i$, $i=1,\ldots, m-1$.
For any $j=1,\ldots,m-1$, one has
$$p_j:=\sum_{i=1}^{j}k_i=\mu_1-\mu_{j+1}\,.
$$
Using $\sum_{k=1}^{m}\mu_k=\frac{m(m-1)}{4}$, one then gets
$$\sum_{j=1}^{m-1}p_j=(m-1)\, \mu_1+\mu_1-\frac{m(m-1)}{4}=m\, \mu_1-\frac{m(m-1)}{4}\,,$$
hence 
$$\mu_1=\frac{m-1}{4}+\frac{1}{m}\sum_{i=1}^{m-1}p_i\,,$$
and
$$\mu_{j+1}=\mu_1-p_j=\frac{m-1}{4}+\frac{1}{m}\left(\sum_{i=1}^{m-1} \,p_i-m\,p_j\right)\,.$$
The expression $\sum_{k=1}^{m} (\mu_{k}^2+k\,\mu_k)$ is a polynomial $F(p_1,\ldots,p_{m-1})$ of the variables 
$p_1$, $p_2\ldots,p_{m-1}$. With the notation $p=(p_1,\ldots,p_{m-1})$, one has
$$\frac{\partial F}{\partial p_i}(p)=-2\mu_{i+1}+m-(i+1)\,.$$
The function $F$ has a unique critical point when $\mu_i=\frac{m-i}{2}$, $i=1,\ldots,m$, which is equivalent to
$p_i=\frac{i}{2}$, or $k_i=\frac{1}{2}$, $i=1,\ldots,m-1$.

\ni Let $p_0=\left(\frac{1}{2},1,\frac{3}{2},\ldots\frac{m-1}{2}\right)$. As
$$\frac{\partial^2 F}{\partial p_i\partial p_j}(p_0)=-\frac{2}{m}\,(1-m\,\delta_{ij})\,,$$
we may apply the result of the previous example. With the same notations,
\begin{align*}
 F(p)-F(p_0)&=\frac{1}{2}\, {}^t (p-p_0) H (p-p_0)\,,\\
 &=\frac{1}{m}\, {}^t (p-p_0)\, Q\begin{pmatrix} 1 &0& 0&0\\
 0&m&0&0\\ &&\ddots&\\0&0&0&m\end{pmatrix}{}^t Q\,(p-p_0)\,,\\
 &=\frac{1}{m}\,\epsilon_1\,\left((m-1)k_1+(m-2)k_2+\cdots +k_{m-1}-\frac{m(m-1)}{4}\right)^2\\
 &\qquad + \epsilon_2\, \left((m-2)k_2+(m-3)k_3\cdots +k_{m-1}-\frac{(m-1)(m-2)}{4}\right)^2\\
 &\qquad +\cdots\\
 &\qquad + \epsilon_{m-i}\, \left(i\,k_{m-i}+\cdots +k_{m-1}-\frac{i(i+1)}{4}\right)^2\\
 &\qquad +\cdots\\
 &\qquad +\epsilon_{m-2}\, \left(2 k_{m-2}+k_{m-1}-\frac{3}{2}\right)^2\\
 &\qquad +\epsilon_{m-1}\, \left(k_{m-1}-\frac{1}{2}\right)^2\,.
\end{align*}
The minimum is obtained if and only if all the squares are minimal. Hence, since the $k_i$ have to be integers, 
\begin{lem}
 The minimum is obtained only when 
 $$(k_1,k_2,\ldots,k_{m-1})=(1,0,1,0,\ldots)\; \text{or}\;
 (0,1,0,1,\ldots)\,.$$
\end{lem}
\proof Note that 
\begin{align*}
 F(p)-F(p_0)&=\frac{1}{4m}\,\epsilon_1\,\left((m-1)(2k_1-1)+(m-2)(2 k_2-1)+\cdots +(2 k_{m-1}-1)
\right)^2\\
 &\qquad +\frac{1}{4}\,\epsilon_2\, \left((m-2)(2k_2-1)+(m-3)(2 k_3-1)\cdots +(2k_{m-1}-1)
 \right)^2\\
 &\qquad +\cdots\\
 &\qquad + \frac{1}{4}\,\epsilon_{m-i}\, \left(i\,(2k_{m-i}-1)+\cdots +(2k_{m-1}-1)\right)^2\\
 &\qquad +\cdots\\
 &\qquad +\frac{1}{4}\,\epsilon_{m-2}\, \left(2 (2k_{m-2}-1)+(2k_{m-1}-1)\right)^2\\
 &\qquad +\frac{1}{4}\,\epsilon_{m-1}\, \left(2k_{m-1}-1\right)^2\,.
\end{align*}
Since the $2 k_i-1$ are odd integers, the result of lemma~\ref{min} shows that the minimum is obtained if and 
only if $(2k_1-1,2k_2-1,\ldots, 2k_{m-1}-1)=(1,-1,1,-1,\ldots)$ or $(-1,1,-1,1,\ldots)$, hence if and only if 
$(k_1,k_2,\ldots,k_{m-1})=(1,0,1,0,\ldots)$ or
 $(0,1,0,1,\ldots)$. \qed

 \ni Hence the $G$-weights $\mu$ for which $\mu_{|\gt_K}=\lambda_{0}$ and $\|\mu+\delta_G\|^2$ is minimal
are given by
\begin{enumerate}
 \item If $m$ is even, $m=2p$,
 \begin{enumerate}
 \item if $(k_1,k_2,\ldots,k_{m-1})=(1,0,1,0,\ldots,1)$, since $p_{2i}=i$, and $p_{2i+1}=i+1$,
 \begin{align*}
  \mu_0:=&\sum_{i=1}^{p} \left(p+\frac{3}{4}-i\right)\, \h{y}_{2i-1}+
  \sum_{i=1}^{p} \left(p-\frac{1}{4}-i\right)\, \h{y}_{2i}\\
  &-\sum_{i=1}^{p} \left(p+\frac{1}{4}-i\right)\, \h{y}_{m+2i-1}-
  \sum_{i=1}^{p} \left(p+\frac{1}{4}-i\right)\, \h{y}_{m+2i}\,.
 \end{align*}
 \item if $(k_1,k_2,\ldots,k_{m-1})=(0,1,0,1,\ldots,0)$, since $p_{2i}=p_{2i+1}=i$,
 \begin{align*}
  \mu_0':=&\sum_{i=1}^{p} \left(p+\frac{1}{4}-i\right)\, \h{y}_{2i-1}+
  \sum_{i=1}^{p} \left(p+\frac{1}{4}-i\right)\, \h{y}_{2i}\\
  &-\sum_{i=1}^{p} \left(p+\frac{3}{4}-i\right)\, \h{y}_{m+2i-1}-
  \sum_{i=1}^{p} \left(p-\frac{1}{4}-i\right)\, \h{y}_{m+2i}\,.
 \end{align*}
Note that $\mu_{0}'=\sigma^{*}(\mu_{0})$ and  $\|\mu_{0}+\delta_G\|^2=\|\mu_{0}'+\delta_G\|^2$.
 \end{enumerate}
 \item If $m$ is odd, $m=2p+1$,
 \begin{enumerate}
 \item if $(k_1,k_2,\ldots,k_{m-1})=(1,0,1,0,\ldots,0)$, since $p_{2i}=i$, and $p_{2i+1}=i+1$,
 \begin{align*}
  \mu_0:=&\sum_{i=1}^{p+1} \left(p+\frac{5}{4}-i-\frac{1}{4m}\right)\, \h{y}_{2i-1}+
  \sum_{i=1}^{p} \left(p+\frac{1}{4}-i-\frac{1}{4m}\right)\, \h{y}_{2i}\\
  &-\sum_{i=1}^{p+1} \left(p+\frac{3}{4}-i+\frac{1}{4m}\right)\, \h{y}_{m+2i-1}-
  \sum_{i=1}^{p} \left(p+\frac{3}{4}-i+\frac{1}{4m}\right)\, \h{y}_{m+2i}\,.
 \end{align*}
 \item if $(k_1,k_2,\ldots,k_{m-1})=(0,1,0,1,\ldots,1)$, since $p_{2i}=p_{2i+1}=i$,
 \begin{align*}
  \mu_0':=&\sum_{i=1}^{p+1} \left(p+\frac{3}{4}-i+\frac{1}{4m}\right)\, \h{y}_{2i-1}+
  \sum_{i=1}^{p} \left(p+\frac{3}{4}-i+\frac{1}{4m}\right)\, \h{y}_{2i}\\
  &-\sum_{i=1}^{p+1} \left(p+\frac{5}{4}-i-\frac{1}{4m}\right)\, \h{y}_{m+2i-1}-
  \sum_{i=1}^{p} \left(p+\frac{1}{4}-i-\frac{1}{4m}\right)\, \h{y}_{m+2i}\,.
\end{align*}
 Note that $\mu_{0}'=\sigma^{*}(\mu_{0})$ and  $\|\mu_{0}+\delta_G\|^2=\|\mu_{0}'+\delta_G\|^2$.
\end{enumerate}
\end{enumerate}
Note that $\mu_0$ and $\mu_0'$ are $G$-dominant, hence we may conclude exactly as in the above case with the 
result of lemma~\ref{low.cas}~:  the square of the first eigenvalue of the Dirac operator is given by 
$$
\frac{1}{4m}\, \langle \mu_{0},\mu_{0}+2\delta_G\rangle+\frac{2m^2-m-1}{16}\,,
$$
hence the result.

\section{The symmetric space $\dsp{\frac{\mathrm{SO}(2p+2q+2)}{\mathrm{SO}(2p+1)\times\mathrm{SO}(2q+1)}}$,
$p\leq q$, $p+q\geq 1$.}
Let $(e_i)_{1\leq i\leq 2p+2q+2}$ be the standard basis of $\bR^{2p+2q+2}$. Let $J$ be the 
diagonal matrix 
$$J=\begin{pmatrix} -I_{2p} &0 &0&0\\
0&I_{2q}&0&0\\
0&0&-1&0\\
0&0&0&1\end{pmatrix}\,.$$
We consider the involution $\sigma$ of the group $\mathrm{SO}(2p+2q+2)$ defined by
$$ A\longmapsto JAJ^{-1}\,.$$
Note that $J$ is orthogonal but $\det(J)=-1$ hence $\sigma$ is not a conjugation in the group.
The connected component of the subgroup of fixed points is isomorphic to
$\mathrm{SO}(2p+1)\times\mathrm{SO}(2q+1)$.

\ni We choose to consider $G=\mathrm{Spin}(2p+2q+2)$ instead. We view $\mathrm{Spin}(2p+1)$ and 
$\mathrm{Spin}(2q+1)$ as subgroups of $G$ by considering
\begin{align*}
 \mathrm{Spin}(2p+1)&= \{v_1\cdots v_{2k}\,;\, v_i\in \mathrm{span}\{e_1,\ldots,e_{2p},e_{2p+2q+1}\}\,;\,
 \|v_i\|=1\}\,,\\
 \intertext{and}
 \mathrm{Spin}(2q+1)&= \{v_1\cdots v_{2\ell}\,;\, v_i\in \mathrm{span}\{e_{2p+1},\ldots,e_{2p+2q},e_{2p+2q+2}\}\,;\,
 \|v_i\|=1\}\,.
\end{align*}

\ni Let $K$ be the connected subgroup of $G$ defined by the image of the morphism 
$$\mathrm{Spin}(2p+1)\times\mathrm{Spin}(2q+1)\longrightarrow \mathrm{Spin}(2p+2q+2)\;;\;
(\varphi,\psi)\longmapsto \varphi\cdot\psi\,.
$$
This group $K$ is the connected component of the subgroup of fixed elements of the (outer) involution
$$\sigma: G\longrightarrow G\;;\; \psi\longmapsto \varphi_{0}\cdot\psi\cdot \varphi_{0}^{-1}\,,$$
where
$$\varphi_{0}:=e_{2p+1}\cdot e_{2p+2}\cdots e_{2p+2q}\cdot e_{2p+2q+2}\,.$$
The decomposition of the Lie algebra $\gspin_{2p+2q+2}$ induced by the involution $\sigma$ is given by
$$\gspin_{2p+2q+2}=(\gspin_{2p+1}\oplus \gspin_{2q+1})\oplus \gp\,,$$
where
$$\gp=\left\{\sum_{\substack{i\in \{1,\ldots,2p,2p+2q+1\}\\
j \in \{2p+1,\ldots,2p+2q,2p+2q+2\}}} \beta_{ij}\, e_{i}\cdot{e_j}\;;\; \beta_{ij}\in \bR\right\}\,.$$
Thus
$$\dim(\gp)=(2p+1)(2q+1)\,.$$
Let $T$ be the standard maximal torus of $G$~:
$$T=\left\{\prod_{k=1}^{p+q+1} \left(\cos\, \beta_k+\sin\, \beta_k\,e_{2k-1}\cdot e_{2k}\right)\;;\; 
\beta_k\in \bR\right\}\,.$$
Then
\begin{align*}
 \gt&=\left\{\sum_{k=1}^{p+q+1} \beta_k\,e_{2k-1}\cdot e_{2k}\;;\; \beta_k\in \bR\right\}\,,\\
 \intertext{and}
 \gt_K&=\left\{\sum_{k=1}^{p} \beta_k\,e_{2k-1}\cdot e_{2k}+
 \sum_{k=p+1}^{p+q} \beta_k\,e_{2k-1}\cdot e_{2k}\;;\; \beta_k\in \bR\right\}\,,\\
 \gt_0&=\mathrm{span}\{e_{2p+2q+1}\cdot e_{2p+2q+2}\}\,.
\end{align*}
Let $\h{x}_k$, $k=1,\ldots, p+q+1$, be the basis of $i\, \gt^{*}$ defined by
$$\h{x}_k(H)=2i\,\beta_k\,,\;\text{for}\; H=\sum_{k=1}^{p+q+1} \beta_k\,e_{2k-1}\cdot e_{2k}\in \gt\,.$$ 
Any element $\mu$ in $i\, \gt^{*}$ of the form 
$$\mu=\sum_{k=1}^{p+q+1} \mu_k\, \h{x}_k\,,\; \mu_k\in \bR\,,$$
is denoted
$$\mu=(\mu_1,\ldots,\mu_{p+q+1})\,.$$
The scalar product on $i\,\gt^{*}$ considered here is given by the scalar product on the Lie algebra $\gspin_{2p+2q+2}$
defined by 
$$\langle X,Y\rangle=-\frac{1}{2}\, \mathrm{Tr}\left(\xi_{*}(X)\xi_{*}(Y)\right)=-\frac{1}{4(p+q)}\,
\mathrm{B}(X,Y)\,,\; X,Y\in \gspin_{2p+2q+2}\,,$$
where $\xi$ is the covering $\mathrm{Spin}(2p+2q+2)\rightarrow \mathrm{SO}(2p+2q+2)$.
For any $\mu=(\mu_1,\ldots,\mu_{p+q+1})$ and any $\mu'=(\mu_1',\ldots,\mu_{p+q+1}')$
$$\langle \mu, \mu'\rangle=\sum_{k=1}^{p+q+1} \mu_{k}\,\mu_{k}'\,.$$
The involution $\sigma^{*}$ of $i\, \gt^{*}$ induced by $\sigma$ is defined by
$$\mu=(\mu_1,\ldots,\mu_{p+q},\mu_{p+q+1})\xrightarrow{\sigma^{*}} 
(\mu_1,\ldots,\mu_{p+q},-\mu_{p+q+1})\,.$$
Note that the vectors $\h{x}'_k:=\h{x}_{k|\gt_K}$, $k=1,\ldots,p+q$, define a basis of $i\, \gt_K^{*}$.
Any $\lambda\in i\, \gt_K^{*}$ of the form $\lambda=\sum_{k=1}^{p+q} \lambda_k\, \h{x}'_k$ is denoted
$$\lambda=(\lambda_1,\ldots,\lambda_{p+q})\,.$$

\subsection{Sets of roots.}
Let 
$$u_k=\frac{1}{2}(e_{2k-1}-i\,e_{2k})\quad \text{and}\quad v_k=\frac{1}{2}(e_{2k-1}+i\,e_{2k})\,,\; k=1,\ldots,
p+q+1\,.$$
The $G$-roots are 
\begin{align*}
 &\h{x}_i+\h{x}_j\,,\quad 1\leq i<j\leq p+q+1\,,\quad \text{with root-vector space $\bC\, u_i\cdot u_j$,}\\
-(&\h{x}_i+\h{x}_j)\,,\quad 1\leq i<j\leq p+q+1\,,\quad \text{with root-vector space $\bC\, v_i\cdot v_j$,}\\
 &\h{x}_i-\h{x}_j\,,\quad 1\leq i<j\leq p+q+1\,,\quad \text{with root-vector space $\bC\, u_i\cdot v_j$,}\\
 -(&\h{x}_i-\h{x}_j)\,,\quad 1\leq i<j\leq p+q+1\,,\quad \text{with root-vector space $\bC\, v_i\cdot u_j$.}\\
\end{align*}
Note that 
\begin{align*}
 \Phi_1&=\left\{\pm(\h{x}_i\pm\h{x}_j)\;;\; 1\leq i<j\leq p\,,\;
p+1\leq i<j\leq p+q\right\}\,,\\
\Phi_2&= \left\{\pm(\h{x}_i\pm\h{x}_j)\;;\; 1\leq i\leq p\,,\;
p+1\leq j\leq p+q\right\}\,,\\
\Phi_3&=\left\{\pm(\h{x}_i\pm\h{x}_{p+q+1})\;;\; 1\leq i\leq p+q\right\}\,.
\end{align*}
We choose positive roots such that
\begin{align*}
 \Phi_{1}^{+}&=\left\{\begin{array}{l}\h{x}_i-\h{x}_j\\\h{x}_i+\h{x}_j\end{array}
\;;\; 1\leq i<j\leq p\,,\;
p+1\leq i<j\leq p+q\right\}\,,\\
\Phi_{2}^{+}&= \left\{\begin{array}{l}\h{x}_i-\h{x}_j\\\h{x}_i+\h{x}_j\end{array}\;;\; 1\leq i\leq p\,,\;
p+1\leq j\leq p+q\right\}\,,\\
\Phi_{3}'^{+}&=\left\{\h{x}_i-\h{x}_{p+q+1}\;;\; 1\leq i\leq p+q\right\}\,,\\
\sigma^{*}(\Phi_{3}'^{+})&=\left\{\h{x}_i+\h{x}_{p+q+1}\;;\; 1\leq i\leq p+q\right\}\,.
\end{align*}
It is easy to see that a system of simple roots is given by
\begin{align*}
\alpha_k&:= \h{x}_k-\h{x}_{k+1}\,,\; k=1,\ldots,p-1\,,\; k=p+1,\ldots, p+q-1\,,(\in \Phi_1^{+})\,,\\
\beta&:=\h{x}_p-\h{x}_{p+1}\,,(\in \Phi_2^{+})\,,\\
\gamma&:=\h{x}_{p+q}-\h{x}_{p+q+1}\,,(\in \Phi_{3}'^{+})\,,\\
\sigma^{*}(\gamma)&:=\h{x}_{p+q}+\h{x}_{p+q+1}\,,(\in \sigma^{*}(\Phi_{3}'^{+}))\,.
\end{align*}
Note that
\begin{align*}
 \langle \alpha_k,\alpha_k\rangle&=\langle \beta,\beta\rangle=\langle \gamma,\gamma\rangle=
 \langle \sigma^{*}(\gamma),\sigma^{*}(\gamma)\rangle=2\,,\\
 &k=1,\ldots,p-1,p+1,\ldots,p+q-1\,,\\
\langle \alpha_k,\alpha_{k+1}\rangle&=\langle \alpha_{p-1},\beta\rangle=\langle \alpha_{p+1},\beta\rangle=
\langle \alpha_{p+q-1},\gamma\rangle=\langle \alpha_{p+q-1},\sigma^{*}\gamma\rangle=-1\,,\\
&k=1,\ldots,p-2,p+1,\ldots,p+q-2\,,
\end{align*}
all the other scalar products being zero. Hence the following (classical) characterization.
\begin{lem} A vector $\mu=(\mu_1,\ldots, \mu_{p+q+1})\in i\, \gt^*$ is $G$-dominant if and only if
$$\mu_1\geq\mu_2\geq \cdots\geq \mu_{p+q}\geq |\mu_{p+q+1}|\,,$$
the $\mu_i$ being all simultaneously integers or half-integers.
\end{lem}
\ni By \eqref{pos.K.roots}, the positive $K$-roots are the restrictions to $\gt_K$ of the positive $G$-roots in
$\Phi_1$ and $\Phi_3'$, hence,
$$\Phi_K^{+}=\left\{\left.\begin{array}{l}\h{x}'_i-\h{x}'_j\\\h{x}'_i+\h{x}'_j \end{array}\;\right\} 
\begin{array}{l}1\leq i<j\leq p\,,\\ p+1\leq i<j\leq p+q\,,\end{array}\;;\; \h{x}'_i\,,\; 1\leq i\leq p+q\right\}\,.
$$
A system of simple roots is given by
$$\theta'_k=\h{x}'_k-\h{x}'_{k+1}\,,\; k=1,\ldots,p-1,p+1,\ldots p+q-1\,,\quad 
\h{x}'_p\quad \text{and}\quad \h{x}'_{p+q}\,.
$$
Note that
\begin{align*}
 \langle \theta'_k,\theta'_k\rangle&=2\,,\; k=1,\ldots,p-1,p+1,\ldots p+q-1\,,\quad 
 \langle \h{x}'_p,\h{x}'_p\rangle=\langle \h{x}'_{p+q},\h{x}'_{p+q}\rangle=1\,,\\
 \langle \theta'_k,\theta'_{k+1}\rangle&=\langle \theta'_{p-1},\h{x}'_p\rangle=
 \langle \theta'_{p+q-1}, \h{x}'_{p+q}\rangle =-1\,,\; k=\begin{cases}1,\ldots,p-2\,,\\
 p+1,\ldots p+q-2\,.\end{cases}
\end{align*}
Hence
\begin{lem} A vector $\lambda=(\lambda_1,\ldots,\lambda_{p+q})\in i\, \gt_K^{*}$ is $K$-dominant
if and only if
$$\lambda_1\geq \lambda_2\geq \cdots \geq \lambda_p\geq 0\,,\quad \text{and}\quad
\lambda_{p+1}\geq \lambda_{p+2}\geq \cdots \geq \lambda_{p+q}\geq 0\,,
$$ 
the $\lambda_i$ for $i=1,\ldots, p$ (resp. for $i=p+1,\ldots, p+q$) being all simultaneously 
integers or half-integers.
\end{lem}
\subsection{A characterization of the highest weights of the spin representation of $K$.}
\label{hw.so}
The explicit determination of the highest weights is far from being simple. Some results
on the decomposition of the spin representation for oriented grassmannians may be found in
 \cite{Kli07}.

\ni By the result of lemma~\ref{h.weight}, any highest weight of the spin representation of $K$ has the form
$$\frac{1}{2}\left(\sum_{\beta\in \Phi_2^{+}} \pm \beta'+\sum_{i=1}^{p+q} \h{x}'_i\right)\,.$$

\ni Now, remark that the $\alpha's$ in $\Phi_1$ and the $\beta's$ in $\Phi_2$ are respectively  compact and
non-compact roots relative to the maximal common torus $T_K$ of the groups $G_1=\mathrm{SO}(2p+2q)$ and $K_1=
\mathrm{SO}(2p)\times \mathrm{SO}(2q)$. Since $G_1/K_1$ is an inner symmetric space, the results of R. Parthasarathy
in \cite{Par} may be applied here.

\ni First the weights of the spin representation of $K_1$ are 
$$\frac{1}{2}\left(\sum_{\beta\in \Phi_2^{+}} \pm \beta'\right)\,,$$
see Remark~2.1 in \cite {Par}. 

\ni Now from lemma~2.2 in \cite {Par} we may conclude
\begin{lem} Any highest weight of the spin representation of $K$ has necessarily the form
$$ w\cdot \delta_{G_1}-\delta_{K_1}+\frac{1}{2}\sum_{i=1}^{p+q} \h{x}'_i\,,\quad w\in W_1\,,$$
where \begin{itemize}
       \item $\delta_{G_1}$ is the half-sum of the
             positive $G_1$-roots, that are the $\alpha's$ in $\Phi_1^{+}$ and the $\beta's$ in $\Phi_2^{+}$,
             and whose set is denoted $\Phi_{G_1}^{+}$,
       \item $\delta_{K_1}$ is the half-sum of the
             positive roots of $K_1$, that are the $\alpha's$ in $\Phi_1^{+}$,
             and whose set is denoted $\Phi_{K_1}^{+}$,
       \item $W_1$ is the subset of the Weyl group $W_{G_1}$ of $G_1$ defined by
       $$W_1=\{w\in W_{G_1}\,;\, w\cdot \Phi_{G_1}^{+}\supset \Phi_{K_1}^{+}\}\,.$$
       \end{itemize}
\end{lem}
\proof Let 
$\lambda=\frac{1}{2}\left(\sum_{\beta\in \Phi_2^{+}} \varepsilon_{\lambda\beta}\, \beta'+\sum_{i=1}^{p+q} \h{x}'_i\right)
$, $\varepsilon_{\lambda\beta}=\pm 1$, be a highest weight of the spin representation of $K$.

\ni For any $\alpha\in \Phi_1$, 
$$\alpha'+\frac{1}{2}\sum_{\beta\in \Phi_2^{+}} \varepsilon_{\lambda\beta}\, \beta'\not =
\frac{1}{2}\sum_{\beta\in \Phi_2^{+}} \pm\, \beta'\,,$$ since otherwise 
$\lambda+\alpha'$ is a weight of the spin representation of $K$, contradicting the fact that $\lambda$
is a highest weight.

\ni So $\frac{1}{2}\sum_{\beta\in \Phi_2^{+}} \varepsilon_{\lambda\beta}\, \beta'$ is a highest weight 
of the spin representation of $K_1$, hence of the form $w\cdot \delta_{G_1}-\delta_{K_1}$, $w\in W_1$, by the 
result of Parthasarathy. \qed

\ni Now let $\lambda$ be a highest weight of that sort. One has 
\begin{align*}
 \delta_{G_1}&=\sum_{i=1}^{p+q}(p+q-i)\, \h{x}_i\,,\\
 \intertext{and}
 \delta_{K_1}&=\sum_{i=1}^p (p-i)\,\h{x}_i+\sum_{i=1}^{p+q} (p+q-i)\, \h{x}_i\,.
\end{align*}
On the other hand, the Weyl group $W_{G_1}$ acts on $i\,\gt_{K}^{*}$ as
$$(\lambda_1,\ldots,\lambda_{p+q})\mapsto (\epsilon_1\, \lambda_{\sigma(1)}, \epsilon_2\, \lambda_{\sigma(2)},
\ldots, \epsilon_{p+q}\, \lambda_{\sigma(p+q)})\,,$$
where $\sigma \in \mathfrak{S}_{p+q}$, $\epsilon_i=\pm 1$, $\epsilon_1\cdots\epsilon_{p+q}=1$, see for instance
\cite{BH3M}. So, for any $w\in W_{G_1}$,
\begin{align*} w\cdot\delta_{G_1}-\delta_{K_1}=
&\sum_{i=1}^{p} \big(\epsilon_{\sigma(i)}\,(p+q-\sigma(i))-(p-i)\big)\, \h{x}_i\\
&+\sum_{i=p+1}^{p+q} \big(\epsilon_{\sigma(i)}\,(p+q-\sigma(i))-(p+q-i)\big)\, \h{x}_i\,,
\end{align*}
where $\sigma\in \mathfrak{S}_{p+q}$. But the dominance conditions verified by $\lambda$ implies
$\epsilon_i=1$, $i=1,\ldots,p-1,p+1,\ldots, p+q-1$, and this is also true for $\epsilon_p$ and
$\epsilon_{p+q}$ since $ \epsilon_p=\epsilon_{p+q}=-1$ implies
$\sigma(p)=p+q=\sigma(p+q)$. Now the dominance conditions also imply that any highest weight $\lambda$ 
of the spin representation of $K$ has necessarily the form
\begin{equation}\label{h.weight.so}
\lambda=\sum_{i=1}^{p} \left(q+i-\sigma(i)+\frac{1}{2}\right)\,\h{x}_i+
\sum_{i=p+1}^{p+q} \left(i-\sigma(i)+\frac{1}{2}\right)\,\h{x}_i\,,
\end{equation}
where $\sigma\in \mathfrak{S}_{p+q}$ verifies $\sigma(i)\leq \sigma(i+1)-1$, $1\leq i\leq p-1$, 
$p+1\leq i\leq p+q-1$. Note that this implies $\sigma(i)\geq i$, $1\leq i\leq p$, and 
$\sigma(i)\leq i$, $p+1\leq i\leq p+q$.

\subsection{The first eigenvalue of the Dirac operator.}

\ni We begin by determining necessary conditions for a $G$-weight $\mu=(\mu_1,\ldots,\mu_{p+q+1})$ in order that
$\mu_{|\gt_K}=(\mu_1,\ldots,\mu_{p+q})$ is a highest weight of the spin representation
and $\|\mu+\delta_G\|^2$ is minimal.

\ni Note first that a vector $\mu=(\mu_1,\ldots,\mu_{p+q+1})\in i\gt^{*}$, is a $G$-weight
if and only if all the $\mu_i$ are all simultaneously integers or half-integers.

\ni So, by \eqref{h.weight.so}, for such a $G$-weight,
$\mu_{|\gt_K}=(\mu_1,\ldots,\mu_{p+q})$ is a highest weight of the spin representation
only if $\mu_{p+q+1}$ is a half-integer. Furthermore, since
$$\delta_G=(p+q,p+q-1,\ldots, 1,0)\,,$$
the condition that $\|\mu+\delta_G\|^2$ is minimal implies then $\mu_{p+q}=\pm \frac{1}{2}$ 
(with the same value in both cases).

\ni So let $\mu$ be such a weight. We first consider the case where $\mu$ is dominant.

\subsubsection{First case: $\mu$ is $G$-dominant.} By \eqref{h.weight.so}, one has
$$\mu+\delta_G=\sum_{i=1}^{p} \left(p+2q+\frac{3}{2}-\sigma(i)\right)\, \h{x}_i
+\sum_{i=p+1}^{p+q} \left(p+q+\frac{3}{2}-\sigma(i)\right)\, \h{x}_i\,,
$$ where $\sigma\in \mathfrak{S}_{p+q}$ verifies $\sigma(i)\leq \sigma(i+1)-1$, $1\leq i\leq p-1$, 
$p+1\leq i\leq p+q-1$, and also, as $\mu$ is $G$-dominant $\sigma(p)\leq \sigma(p+1)+q-1$.

\ni So
\begin{align*}
\|\mu+\delta_G\|^2 &= \sum_{i=1}^{p} \left(p+2q+\frac{3}{2}-\sigma(i)\right)^2+
\sum_{i=p+1}^{p+q} \left(p+q+\frac{3}{2}-\sigma(i)\right)^2\,,\\
&=\sum_{i=1}^{p+q} \sigma(i)^2 -2\,\left(p+q+\frac{3}{2}\right)\,\sum_{i=1}^{p+q} \sigma(i)
-2\,q\,\sum_{i=1}^{p} \sigma(i)\\
& \qquad\qquad +\sum_{i=1}^{p} \left(p+2q+\frac{3}{2}\right)^2+
\sum_{i=p+1}^{p+q} \left(p+q+\frac{3}{2}\right)^2\,,\\
&=\sum_{i=1}^{p+q} i^2 -2\,\left(p+q+\frac{3}{2}\right)\,\sum_{i=1}^{p+q} i
-2\,q\,\sum_{i=1}^{p} \sigma(i)\\
& \qquad\qquad +\sum_{i=1}^{p} \left(p+2q+\frac{3}{2}\right)^2+
\sum_{i=p+1}^{p+q} \left(p+q+\frac{3}{2}\right)^2\,.
\end{align*}
Hence 
$$ \|\mu+\delta_G\|^2\;\text{is minimal} \Longleftrightarrow \sum_{i=1}^{p} \sigma(i)\; \text{is maximal}\,.$$
Note that the conditions $\sigma(i)\leq \sigma(i+1)-1$, $1\leq i\leq p-1$ imply that $\sum_{i=1}^{p} \sigma(i)$ is 
maximal if and only if $\sigma(p)$ is maximal and $\sigma(p-1)=\sigma(p)-1$,..., $\sigma(1)=\sigma(p)-p+1$.
But $\sigma(p)$ can not be $>q$ in that case. Indeed, if $\sigma(p)>q$, then, as we suppose $p\leq q$, 
$\sigma(1)>1$. So $\sigma(p+1),\ldots \sigma(p+q)$ belong to the set 
$\{1,\ldots,\sigma(1)-1,\sigma(p)+1,\ldots,p+q\}$, and then the conditions
$\sigma(i)\leq \sigma(i+1)-1$, $p+1\leq i\leq p+q-1$, imply $\sigma(p+1)=1$. But that contradicts 
the condition $\sigma(p)\leq \sigma(p+1)+q-1$.

\ni Hence $\|\mu+\delta_G\|^2$ is minimal if and only if
$$\begin{array}{ccc} \sigma(1)&=&q-p+1\,,\\
   \sigma(2)&=&q-p+2\,,\\
   \vdots&\vdots&\vdots\\
   \sigma(p)&=&q\,,\\
   \sigma(p+1)&=&1\,,\\
   \sigma(p+2)&=&2\,,\\
   \vdots&\vdots&\vdots\\
   \sigma(q)&=&q-p\,,\\
   \sigma(q+1)&=&q+1\,,\\
   \vdots&\vdots&\vdots\\
   \sigma(p+q)&=&p+q\,,
  \end{array}$$
So we may conclude 
\begin{lem} If $\mu$ is a $G$-dominant weight such that $\mu_{|\gt_K}$ is a 
highest weight $\lambda$ of the spin representation and $\|\mu+\delta_G\|^2$ is minimal,
then necessarily,
$$\mu=(\underbrace{p+\frac{1}{2},\ldots,p+\frac{1}{2}}_{q},
\underbrace{\frac{1}{2},\ldots,\frac{1}{2}}_{p},\pm\frac{1}{2})\,,$$
and 
$$\lambda=(\underbrace{p+\frac{1}{2},\ldots,p+\frac{1}{2}}_{q},
\underbrace{\frac{1}{2},\ldots,\frac{1}{2}}_{p})\,.$$
\end{lem}
\subsubsection{Second case: $\mu$ is not $G$-dominant.}
In that case, one obtains also that
$$ \|\mu+\delta_G\|^2\;\text{is minimal} \Longleftrightarrow \sum_{i=1}^{p} \sigma(i)\; \text{is maximal}\,,$$
but now $\sigma\in \mathfrak{S}_{p+q}$ has only to verify the conditions
$\sigma(i)\leq \sigma(i+1)-1$, $1\leq i\leq p-1$, 
$p+1\leq i\leq p+q-1$. In that case $\sigma(p)$ is maximal if and only if $\sigma(p)=p+q$, 
hence $\|\mu+\delta_G\|^2$ is minimal if and only if
$$\begin{array}{ccc} \sigma(1)&=&q+1\,,\\
   \sigma(2)&=&q+2\,,\\
   \vdots&\vdots&\vdots\\
   \sigma(p)&=&p+q\,,\\
   \sigma(p+1)&=&1\,,\\
   \sigma(p+2)&=&2\,,\\
   \vdots&\vdots&\vdots\\
   \sigma(p+q)&=&q\,,
  \end{array}$$
hence 
$$\mu=(\underbrace{\frac{1}{2},\ldots,\frac{1}{2}}_{p},
\underbrace{p+\frac{1}{2},\ldots,p+\frac{1}{2}}_{q},\pm\frac{1}{2})\,,$$
Now, as $\mu$ is conjugate under $W_G$ to one and only one $G$-dominant weight, one sees that since
the Weyl group $W_G$ acts on $i\,\gt^{*}$ as
$$(\mu_1,\ldots,\mu_{p+q+1})\mapsto (\epsilon_1\, \mu_{\sigma(1)}, \epsilon_2\, \mu_{\sigma(2)},
\ldots, \epsilon_{p+q+1}\, \mu_{\sigma(p+q+1)})\,,$$
where $\sigma \in \mathfrak{S}_{p+q+1}$, $\epsilon_i=\pm 1$, $\epsilon_1\cdots\epsilon_{p+q+1}=1$,
$\mu$ is conjugate under the Weyl group to the $G$-dominant weight met in the above case
$$(\underbrace{p+\frac{1}{2},\ldots,p+\frac{1}{2}}_{q},
\underbrace{\frac{1}{2},\ldots,\frac{1}{2}}_{p},\pm\frac{1}{2})\,.$$
Finally, we may conclude that
\begin{lem}\label{dom.w.so} If a $G$-dominant weight $\mu$ verifies the spin condition and is such that
$\|\mu+\delta_G\|^2$ is minimal,
then necessarily,
$$\mu=(\underbrace{p+\frac{1}{2},\ldots,p+\frac{1}{2}}_{q},
\underbrace{\frac{1}{2},\ldots,\frac{1}{2}}_{p},\pm\frac{1}{2})\,.$$
\end{lem}
\ni In order to conclude, we only have to verify that 
\begin{lem} The vector $\lambda\in i\, \gt_K^{*}$ defined by
 $$\lambda=(\underbrace{p+\frac{1}{2},\ldots,p+\frac{1}{2}}_{q},
\underbrace{\frac{1}{2},\ldots,\frac{1}{2}}_{p})\,,$$
is a highest weight of the spin representation of $K$.
\end{lem}
\proof As $\gp$ is odd-dimensional, we have to use the description of spinors given in 
\S~\ref{odd.case}. By the choice of an orthonormal basis of $\gp$ such as \eqref{onbp}, the group 
$\mathrm{SO}(\gp)$ is identified with $\mathrm{SO}(4pq+2(p+q)+1)$, which is itself embedded in the group
$\mathrm{SO}(4pq+2(p+q)+2)$ in such a way that it acts trivially on the last vector of the standard 
basis of $\bR^{4pq+2(p+q)+2}$.

\ni  Let $(Z,\ol{Z})$ be the Witt basis defined by the two orthonormal vectors
on which $\mathrm{SO}(\gp)$ acts trivially. Using the notations of \S~\ref{w.spin}, 
a basis of spinors is given by considering the vectors
\begin{align*}
 &E_I\cdot V_J\cdot \ol{w}\,,\quad \text{if $\# I+\#J$ is even,}\\
 &E_I\cdot V_J\cdot Z\cdot \ol{w}\,,\quad \text{if $\# I+\#J$ is odd.}
\end{align*}
Denoting for short $E_{ij}=E_{\h{x}_i+\h{x}_j}$, $E'_{ij}=E_{\h{x}_i-\h{x}_j}$, $1\leq i\leq p$, $p+1\leq j\leq p+q$,
and $V_j=E_{\h{x}_j-\h{x}_{p+q+1}}$, $1\leq j\leq p+q$, one gets from \eqref{weight}, that the spinor
\begin{align*}v_\lambda&:=\prod_{\substack{{1\leq i\leq p}\\{p+1\leq j\leq p+q}}} E_{ij}\cdot 
\prod_{\substack{{1\leq i\leq p}\\{q+1\leq j\leq p+q}}} E'_{ij}\cdot V_1\cdots V_{p+q}\cdot \ol{w}\,,
\; \text{if $(p+1)(p+q)$ is even,}\\
&:=\prod_{\substack{{1\leq i\leq p}\\{p+1\leq j\leq p+q}}} E_{ij}\cdot 
\prod_{\substack{{1\leq i\leq p}\\{q+1\leq j\leq p+q}}} E'_{ij}\cdot V_1\cdots V_{p+q}\cdot Z\cdot \ol{w}\,,
\; \text{if $(p+1)(p+q)$ is odd,}
\end{align*}
is a weight-vector for the weight $\lambda$.

\ni Now, it may be checked that the root-vector $u_k\cdot v_{k+1}$, associated to the simple root $\theta'_k=\h{x}'_k-\h{x}'_{k+1}$,
acts on spinors by a linear combination of 
\begin{align*}
 &\sum_{j=p+1}^{p+q}\ol{E_{k+1\,j}}\cdot E_{kj}\,,\;\sum_{j=p+1}^{p+q}\ol{E'_{k+1\,j}}\cdot E'_{kj}\,,\quad 
 \text{and}\quad \ol{V_{k+1}}\cdot V_k\,,\quad \text{if $1\leq k\leq p-1$,}\\
 &\sum_{i=1}^{p}\ol{E_{i k+1}}\cdot E_{ik}\,,\;\sum_{i=1}^{p}\ol{E'_{ik}}\cdot E'_{i k+1}\,,\quad 
 \text{and}\quad
  \ol{V_{k+1}}\cdot V_k\,,\quad \text{if $p+1\leq k\leq p+q-1$.}\\
\end{align*}
So $v_\lambda$ is killed by the action of $u_k\cdot v_{k+1}$ since all the $E_{ij}$'s and all the $V_j$'s occur
in the expression of $v_\lambda$, and either $E'_{kj}$ (resp. $E'_{i k+1}$) occur in that expression or
both $E'_{kj}$ and $E'_{k+1\,j}$,
(resp. $E'_{i k+1}$ and $E'_{ik}$) do not occur in the expression.

\ni In the same way, $v_\lambda$ is killed by the action of the root-vector $u_p\cdot e_{2p+2q+1}$,
associated to the simple root $\h{x}'_p$, since that root-vector 
 acts on spinors by a linear combinations of
$$\sum_{j=p+1}^{p+q} E_{pj}\cdot \ol{V_j}\,,\; \sum_{j=p+1}^{p+q} E'_{pj}\cdot V_j\,,\quad\text{and}
\quad V_p\cdot(e_{2p+2q+1}\cdot e_{2p+2q+2})\,.$$

\ni Finally, $v_\lambda$ is also killed by the action of the root-vector $u_{p+q+1}\cdot e_{2p+2q+2}$,
associated to the simple root $\h{x}'_{p+q}$, since that root-vector 
 acts on spinors by a linear combinations of
 $$\sum_{i=1}^{p} E_{i\, p+q}\cdot \ol{V_i}\,,\; \sum_{i=1}^{p} E'_{i\, p+q}\cdot V_i\,,\quad\text{and}
\quad V_{p+q}\cdot(e_{2p+2q+1}\cdot e_{2p+2q+2})\,.$$
Thus $v_\lambda$ is a maximal vector, and so $\lambda$ is a highest weight of the spin representation.\qed

\ni  So $\mu_{\pm}=\displaystyle{(\underbrace{p+\frac{1}{2},\ldots,p+\frac{1}{2}}_{q},
\underbrace{\frac{1}{2},\ldots,\frac{1}{2}}_{p},\pm\frac{1}{2})}$ are two $G$-dominant weights verifying the spin
condition and such that $\|\mu+\delta_G\|^2$ is minimal among all the $G$-dominant weights $\mu$
verifying the spin condition.

\ni Hence the square of the first eigenvalue $\lambda$ of the Dirac operator is given by
\begin{align*}
 \lambda^2&=\frac{1}{4(p+q)}\, \langle \mu_{+},\mu_{+}+2\delta_G\rangle+ \frac{(2p+1)(2q+1)}{16}\,,\\
 &=\frac{1}{16\,(p+q)}\, \big(8pq\,(2p+q+1)+4p\,(p+1)+4q\,(q+1)+1\big)\,.
\end{align*}

\ni Note that for $p=0$, we retrieve\footnote{in that case, $\Phi_2=\emptyset$, so the spin representation 
of $K$ has only one highest weight.} the value of the square of the first eigenvalue of the Dirac operator on the 
standard sphere $S^{2q+1}$~:
$$\lambda^2=\frac{1}{4q}\, \left(\frac{2q+1}{2}\right)^2\,,$$
cf. \cite{Sul}.

\section{The symmetric space $\dsp{\frac{\mathrm{E}_6}{\mathrm{F}_4}}$.}

\ni We use here the results of Murakami \cite{Mur65} (for an outline, see for instance the section ``{non-inner
involutions}'' in chapter~3 of \cite{BR90}).

\ni Explicit computations can be made (and are obtained)
with the help of the programs GAP3, \cite{Gap3}, and LiE, \cite{LiE}.

\ni First, the existence of non-inner involutions corresponds with the existence of non-trivial
symmetries of the Dynkin diagram.

\ni The Dynkin diagram of $\mathrm{E}_6$ is
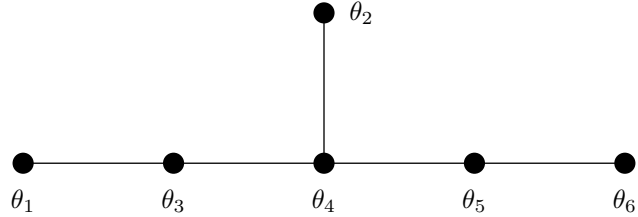
\begin{center}
 \begin{tikzpicture}
  \draw[line width=0.5pt] (-4,0)--(4,0);\draw[line width=0.5pt] (0,0)--(0,2);
  \fill (-4,0) circle(4pt); \draw (-4,-0.5) node{$\theta_1$};
  \fill (-2,0) circle(4pt); \draw (-2,-0.5) node{$\theta_3$};
  \fill (0,0) circle(4pt); \draw (0,-0.5) node{$\theta_4$};
  \fill (2,0) circle(4pt); \draw (2,-0.5) node{$\theta_5$};
  \fill (4,0) circle(4pt); \draw (4,-0.5) node{$\theta_6$};
  \fill (0,2) circle(4pt); \draw (0.5,2) node{$\theta_2$};
  \end{tikzpicture}
\end{center}
There is only one non-trivial symmetry $s$ given by
$$\begin{array}{ccc}
   s(\theta_1)=\theta_6\,, &s(\theta_2)=\theta_2\,, &s(\theta_3)=\theta_5\,,\\
   s(\theta_4)=\theta_4\,, &s(\theta_5)=\theta_3\,, &s(\theta_6)=\theta_1\,.
  \end{array}
$$
The symmetry $s$ is extended by linearity to an involution $\sigma^{*}$ of $i\, \gt^{*}$, which itself 
induces an involution $\sigma_{*}$ of $\gt$, by means of the scalar product (re-normalized here in such a way 
that all simple roots $\theta$ verify $\|\theta\|^2=2$).

\ni Now, choosing a root-vector $E_\theta$ for each simple root $\theta$, $\sigma_{*}$ is extended to the span of 
these vectors by 
\begin{equation}\label{def.sigma}
 \sigma_*(E_\theta)=E_{\sigma_{*}(\theta)}\,.
\end{equation}

\ni Finally $\sigma_{*}$ can be uniquely extended to a non-inner involution of $\gg$ (cf. \S~14.2 in \cite{Hum}).

\ni A first outer symmetric space structure is obtained by considering the connected subgroup $K$ of $\mathrm{E}_6$
whose Lie algebra is the set $\{X\in \gg\,;\, \sigma_{*}(X)=X\}$. It is a simple group, and a system of simple roots
(relatively to the maximal torus $T_K=K\cap T$) is obtained by considering the restriction to $\gt_K$ of the simple roots 
$\theta_i$, $1\leq i\leq 6$, (cf. Proposition~3.20 in \cite{BR90}).

\ni The group $\mathrm{E}_6$ has $36$ positive roots. The positive roots $\theta$ such that
$\sigma^{*}(\theta)=\theta$ all belong to $\Phi_1^{+}$ by \eqref{def.sigma}. The partition of the set of positive
roots is given by
\begin{align*}
 \Phi_1^{+}&=\{\theta_2,\theta_4,\theta_2+\theta_4,\theta_3+\theta_4+\theta_5, \theta_2+\theta_3+
 \theta_4+\theta_5,\theta_1+\theta_3+\theta_4+\theta_5+\theta_6\,,\\
 &\qquad \theta_2+\theta_3+2\,\theta_4+\theta_5,\theta_1+\theta_2+\theta_3+\theta_4+\theta_5+\theta_6\,,\\
 &\qquad \theta_1+\theta_2+\theta_3+2\,\theta_4+\theta_5+\theta_6,
 \theta_1+\theta_2+2\,\theta_3+2\,\theta_4+2\,\theta_5+\theta_6\,,\\
 &\qquad \theta_1+\theta_2+2\,\theta_3+3\,\theta_4+2\,\theta_5+\theta_6,
 \theta_1+2\,\theta_2+2\,\theta_3+3\,\theta_4+2\,\theta_5+\theta_6\}\,,\\
\Phi_2^{+}&=\emptyset\,,\\
\Phi_3'^{+}&=\{\theta_1,\theta_3,\theta_1+\theta_3,\theta_3+\theta_4,\theta_1+\theta_3+\theta_4,
\theta_2+\theta_3+\theta_4,\theta_1+\theta_2+\theta_3+\theta_4\,,\\
&\qquad \theta_1+\theta_3+\theta_4+\theta_5,\theta_1+\theta_2+\theta_3+\theta_4+\theta_5,
\theta_1+\theta_2+\theta_3+2\,\theta_4+\theta_5\,,\\
&\qquad \theta_1+\theta_2+2\,\theta_3+2\,\theta_4+\theta_5,
\theta_1+\theta_2+2\,\theta_3+2\,\theta_4+\theta_5+\theta_6\}\,,\\
\sigma^{*}(\Phi_3'^{+})&=\{\theta_6,\theta_5,\theta_5+\theta_6,\theta_4+\theta_5,\theta_4+\theta_5+\theta_6,
\theta_2+\theta_4+\theta_5,\theta_2+\theta_4+\theta_5+\theta_6\,,\\
&\qquad \theta_3+\theta_4+\theta_5+\theta_6,\theta_2+\theta_3+\theta_4+\theta_5+\theta_6,
\theta_2+\theta_3+2\,\theta_4+\theta_5+\theta_6\,,\\
&\qquad \theta_2+\theta_3+2\,\theta_4+2\,\theta_5+\theta_6,
\theta_1+\theta_2+\theta_3+2\,\theta_4+2\,\theta_5+\theta_6\}\,.
\end{align*}
The set of positive $K$-roots is
$$\Phi_K^{+}=\{\theta_{|\gt_K};\, \theta \in \Phi_1^{+}\cup \Phi_3'^{+}\}\,,$$
and a system of simple $K$-roots is given by
$$\theta'_1=\theta'_6=\frac{1}{2}(\theta'_1+\theta'_6)\,,\quad
\theta'_2\,,\quad \theta'_3=\theta'_5=\frac{1}{2}(\theta'_3+\theta'_5)\,,\quad
\theta'_4\,.$$
Note that 
\begin{equation}
\begin{split}\dim \gg=78\,,\quad \dim\gk=52\,,\quad \dim \gp=26\,,\nonumber\\
 \quad \dim \gt=6\,,\quad \dim \gt_K=4\,,\quad\dim \gt_0=2\,.  
  \end{split}
\end{equation}
One has
\begin{equation}
 \begin{split}
  \|\theta_1'\|^2=\|\theta_3'\|^2=1\,,\quad \|\theta_2'\|^2=\|\theta_4'\|^2=2\,,\nonumber\\
  \langle \theta_1',\theta_3'\rangle=-\frac{1}{2}\,,\quad \langle \theta_2',\theta_4'\rangle=
  \langle \theta_3',\theta_4'\rangle=-1\,,
 \end{split}
\end{equation}
all the other scalar products being zero.
Hence the Dynkin diagram of $K$ is
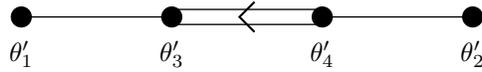
\begin{center}
 \begin{tikzpicture}[line width=0.5pt]
  \draw (-3,0)--(-1,0);\draw (1,0)--(3,0);
  \draw (-1,0.1)--(1,0.1);\draw (-1,-0.1)--(1,-0.1);
  \draw[thick] (0.1,0.2)--(-0.1,0)--(0.1,-0.2);
  \fill (-3,0) circle(4pt); \draw (-3,-0.5) node{$\theta_1'$};
  \fill (-1,0) circle(4pt); \draw (-1,-0.5) node{$\theta_3'$};
  \fill (1,0) circle(4pt); \draw (1,-0.5) node{$\theta_4'$};
  \fill (3,0) circle(4pt); \draw (3,-0.5) node{$\theta_2'$};
  \end{tikzpicture}
\end{center}
Setting
$$\alpha_1=\theta'_2\,,\quad \alpha_2=\theta'_4\,,\quad \alpha_3=\theta'_3\,,
\quad\text{and}\quad \alpha_4=\theta'_1\,,$$
this is the ``{classical}'' Dynkin diagram of the group $\mathrm{F_4}$.
In that basis, the Cartan matrix is 
$$\begin{pmatrix} 2&-1&0&0\\-1&2&-1&0\\0&-2&2&-1\\0&0&-1&2
\end{pmatrix}\,.
$$
In the following, $\omega_1,\omega_2,\omega_3,\omega_4$ are the fundamental weights
associated with $\alpha_1,\alpha_2$, $\alpha_3,\alpha_4$.

\subsection{Highest weights of the spin representation of $K$.}
Since $\Phi_2=\emptyset$, and $\dim \gt_0=2$, we may conclude from lemma~\ref{h.weight},
\begin{lem} The spin representation of $K$ has only one highest weight $\lambda$ 
with multiplicity $2$,
 $$\lambda=5\, \theta_1'+3\, \theta_2'+9\, \theta_3'+
 6\,\theta_4'\,.$$
\end{lem}
\proof By the result of lemma~\ref{h.weight}, 
\begin{align*}
 \lambda=\frac{1}{2}\sum_{\gamma\in \Phi_3'^{+}} \gamma'&=5\, \theta_1'+3\, \theta_2'+9\, \theta_3'+
 6\,\theta_4'\,,\\
 &=3\,\alpha_1+6\, \alpha_2+9\, \alpha_3+5\,\alpha_4=\omega_3+\omega_4\,.
\end{align*}
Any irreducible module with highest weight $\omega_3+\omega_4$ has dimension $2^{12}$.
Since $\dim \gt_0=2$, one knows from \S~\ref{w.spin} that the multiplicity of the weight is at least $2$. 
Now $2\times 2^{12}=2^{13}=2^{\dim(\gp)/2}=\dim \Sigma$, hence the result.\qed

\subsection{The first eigenvalue of the Dirac operator.}

\ni As we did before, we first determine $G$-weights $\mu$ (non necessarily dominant) such that
$\mu_{|\gt_K}=\lambda$ and $\|\mu+\delta_G\|^2$ is minimal.

\ni Let $\mu=\sum_{i=1}^6 \mu_i\, \theta_i\in i\, \gt^{*}$. First
$$\mu_{|\gt_K}=\lambda\Longleftrightarrow \begin{cases}\mu_1+\mu_6&=5\,,\\
\mu_2&=3\,,\\ \mu_3+\mu_5&=9\,,\\ \mu_4&=6\,.\end{cases}$$
As $\|\theta_i\|^2=2$, $1\leq i\leq 6$,
$\mu$ is a $G$-weight (resp. dominant $G$-weight) if and only if $\langle \mu,\theta_i\rangle\in \bZ$, (resp. 
$\bN$), hence $\mu_{|\gt_K}=\lambda$ and $\mu$ is a $G$-weight (resp. dominant $G$-weight) if and 
only if 
$$ \begin{cases}\mu_1+\mu_6&=5\,,\\
\mu_2&=3\,,\\ \mu_3+\mu_5&=9\,,\\ \mu_4&=6\,,\end{cases}
\quad \text{and}\quad
\left.\begin{array}{r}
2\,\mu_1-\mu_3\\
-\mu_1+2\,\mu_3-6\\
\mu_1-2\,\mu_3+7\\
-2\,\mu_1+\mu_3+1
\end{array}\right\}\in \bZ\;\text{(resp. $\bN$).}
$$
Setting
$$\begin{cases}k:=2\,\mu_1-\mu_3\,, \\
l:=-\mu_1+2\, \mu_3-6\,,\end{cases}
\left(\Longleftrightarrow
\begin{cases}\mu_1=\frac{2k+l+6}{3}\,,\\
\mu_2=\frac{k+2l+12}{3}\,,\end{cases}\right)$$
the last condition is equivalent to
$$\text{$k$ and $l\in \bZ$}\;\quad \text{(resp. $k=0$ or $1$ and $l=0$ or $1$).}$$
Viewing then $\|\mu+\delta_G\|^2$ as a polynomial $f(k,l)$ of the variables $k$ and $l$,
one gets 
$$\frac{\partial f}{\partial k}(k,l)=2\,(2\mu_1-5)\quad \text{and}\quad
\frac{\partial f}{\partial l}(k,l)=2\,(2\mu_3-9)\,,$$
Thus $f(k,l)$ has only one critical point $(\frac{1}{2},\frac{1}{2})$. Now
$$\frac{\partial^2 f}{\partial k^2}(k,l)=\frac{8}{3}\,,\quad 
\frac{\partial^2 f}{\partial k\,\partial l}(k,l)=\frac{4}{3}\quad\text{and}
\quad \frac{\partial^2 f}{\partial l^2}(k,l)=\frac{8}{3}\,,$$
So
$$f(k,l)-f\left(\frac{1}{2},\frac{1}{2}\right)=\frac{2}{3}\,\left((k-\frac{1}{2})^2+(l-\frac{1}{2})^2+(k+l-1)^2\right)\,.$$
Hence among the $G$-weights $\mu$ such that $\mu_{|\gt_K}=\lambda$, the minimum of $\|\mu+\delta_G\|^2$ is obtained
if and only if 
\begin{align*}
 &(k,l)=(0,1)\,, \left(\Leftrightarrow (\mu_1,\mu_3)=\left(\frac{7}{3},\frac{14}{3}\right)\right)\,,\\
 \intertext{or}
 &(k,l)=(1,0)\,, \left(\Leftrightarrow (\mu_1,\mu_3)=\left(\frac{8}{3},\frac{13}{3}\right)\right)\,.
\end{align*}
By the above remarks, the corresponding two weights 
\begin{align*}
\mu_1&=\frac{7}{3}\,\theta_1+3\, \theta_2+\frac{14}{3}\, \theta_3+6\,\theta_4+\frac{13}{3}\,
\theta_5+\frac{8}{3}\,\theta_6\,,\\
\intertext{and}
\mu_2&=\frac{8}{3}\,\theta_1+3\, \theta_2+\frac{13}{3}\, \theta_3+6\,\theta_4+\frac{14}{3}\,
\theta_5+\frac{7}{3}\,\theta_6\,,
\end{align*}
are $G$-dominant. Hence we may conclude exactly as we did it in the first example with the result of 
lemma~\ref{low.cas}. We have first to note that, for the scalar product $\langle\,,\,\rangle_K$ induced
by the Killing form sign-changed, the ``{strange formula}'' of Freudenthal and de Vries, \cite{FdV}, gives
$$\|\delta_G\|_K^2=\frac{\dim(\gg)}{24}=\frac{13}{4}\,,$$
whereas for our choice of scalar product $\|\delta_G\|^2=78$, hence $\langle\,,\,\rangle_K=\frac{1}{24}\,
\langle\,,\,\rangle$.

\ni So, finally, the square of the first eigenvalue of the Dirac operator is given by
$$\frac{1}{24} \langle \mu_1+2\,\delta_G,\mu_1\rangle +\frac{\dim(\gp)}{16}=\frac{20}{9}+\frac{13}{8}=
\frac{277}{72}\,.$$

\section{The symmetric space $\displaystyle{\frac{\mathrm{E}_6}{\mathrm{Sp}_4}}$.}
By the result of Murakami, there also exists a complementary non-inner involution on the Lie algebra of 
$\mathrm{E}_6$, which is not conjugate to the above non-inner involution $\sigma$ and  is defined on the Lie algebra 
by 
$$\sigma_{*}'=\sigma_{*}\circ \mathrm{Ad}_{\exp(\pi\, \xi_2)}\,,$$
where $\xi_2\in \gt$ is defined by $\theta_j(\xi_2)=i\, \delta_{j2}$, cf. Theorem~3.25 in \cite{BR90}.

\ni A new outer symmetric space structure is obtained by considering the connected subgroup $K'$ of $\mathrm{E}_6$
whose Lie algebra is the set $\{X\in \gg\,;\, \sigma_{*}'(X)=X\}$.

\ni Note that as $\xi_2\in \gt$, $\sigma'^{*}=\sigma^{*}$.

\ni The partition of the set of positive roots is now given by
\begin{align*}
\Phi_1^{+}&=\{\theta_4,\theta_3+\theta_4+\theta_5,\theta_1+\theta_3+\theta_4+\theta_5+\theta_6,\\
&\qquad 
\theta_1+2\,\theta_2+2\,\theta_3+3\,\theta_4+2\,\theta_5+\theta_6\}\,,\\
\Phi_2^{+}&=\{\theta_2,\theta_2+\theta_4,\theta_2+\theta_3+\theta_4+\theta_5,
\theta_2+\theta_3+2\,\theta_4+\theta_5\,,\\
&\qquad \theta_1+\theta_2+\theta_3+\theta_4+\theta_5+\theta_6,
\theta_1+\theta_2+\theta_3+2\,\theta_4+\theta_5+\theta_6,,\\
& \qquad \theta_1+\theta_2+2\,\theta_3+2\,\theta_4+2\,\theta_5+\theta_6,
\theta_1+\theta_2+2\,\theta_3+3\,\theta_4+2\,\theta_5+\theta_6\}\,,\\
\Phi_3'^{+}&=\{\theta_1,\theta_3,\theta_1+\theta_3,\theta_3+\theta_4,\theta_1+\theta_3+\theta_4,
\theta_2+\theta_3+\theta_4,\theta_1+\theta_2+\theta_3+\theta_4\,,\\
&\qquad \theta_1+\theta_3+\theta_4+\theta_5,\theta_1+\theta_2+\theta_3+\theta_4+\theta_5,
\theta_1+\theta_2+\theta_3+2\,\theta_4+\theta_5\,,\\
&\qquad \theta_1+\theta_2+2\,\theta_3+2\,\theta_4+\theta_5,
\theta_1+\theta_2+2\,\theta_3+2\,\theta_4+\theta_5+\theta_6\}\,,\\
\sigma^{*}(\Phi_3'^{+})&=\{\theta_6,\theta_5,\theta_5+\theta_6,\theta_4+\theta_5,\theta_4+\theta_5+\theta_6,
\theta_2+\theta_4+\theta_5,\theta_2+\theta_4+\theta_5+\theta_6\,,\\
&\qquad \theta_3+\theta_4+\theta_5+\theta_6,\theta_2+\theta_3+\theta_4+\theta_5+\theta_6,
\theta_2+\theta_3+2\,\theta_4+\theta_5+\theta_6\,,\\
&\qquad \theta_2+\theta_3+2\,\theta_4+2\,\theta_5+\theta_6,
\theta_1+\theta_2+\theta_3+2\,\theta_4+2\,\theta_5+\theta_6\}\,.
\end{align*}
The set of positive $K'$-roots is
$$\Phi_K'^{+}=\{\theta_{|\gt_K'};\, \theta \in \Phi_1^{+}\cup \Phi_3'^{+}\}\,,$$
and a system of simple $K'$-roots is given by
$$\alpha_1':=\theta_2'+\theta_4'+\frac{1}{2}(\theta_3'+\theta_5')\,,\; \alpha_2':=\frac{1}{2}(\theta_1'+
\theta_6')\,,\; \alpha_3':=\frac{1}{2}(\theta_3'+\theta_5')\,,\;
\alpha_4'=\theta_4'\,,$$
(cf. Theorem~3.25 in \cite{BR90}).

\ni Note that here
\begin{equation}
\begin{split}\dim \gg=78\,,\quad \dim\gk'=36\,,\quad \dim \gp=42\,,\nonumber\\
 \quad \dim \gt=6\,,\quad \dim \gt_{K'}=4\,,\quad\dim \gt_0=2\,.  
  \end{split}
\end{equation}
One has
\begin{equation}
 \begin{split}
  \|\alpha_1'\|^2=\|\alpha_2'\|^2=\|\alpha_3'\|^2=1\,,\quad \|\alpha_4'\|^2=2\,,\nonumber\\
   \langle \alpha_1',\alpha_2'\rangle=
  \langle \alpha_2',\alpha_3'\rangle=-\frac{1}{2}\,,\quad 
  \langle \alpha_3',\alpha_4'\rangle=-1\,,
 \end{split}
\end{equation}
all the other scalar products being zero. The Cartan matrix is
$$\begin{pmatrix} 2&-1&0&0\\-1&2&-1&0\\0&-1&2&-1\\0&0&-2&2\end{pmatrix}\,,$$
so the Dynkin diagram of $K'$ is 
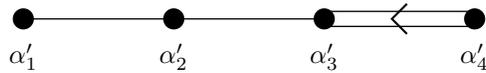
\begin{center}
 \begin{tikzpicture}[line width=0.5pt]
  \draw (-3,0)--(1,0);
  \draw (1,0.1)--(3,0.1);\draw (1,-0.1)--(3,-0.1);
  \draw[thick] (2.1,0.2)--(1.9,0)--(2.1,-0.2);
  \fill (-3,0) circle(4pt); \draw (-3,-0.5) node{$\alpha_1'$};
  \fill (-1,0) circle(4pt); \draw (-1,-0.5) node{$\alpha_2'$};
  \fill (1,0) circle(4pt); \draw (1,-0.5) node{$\alpha_3'$};
  \fill (3,0) circle(4pt); \draw (3,-0.5) node{$\alpha_4'$};
  \end{tikzpicture}
\end{center}
This is the Dynkin diagram of $\mathrm{Sp}_4$.

\subsection{Highest weights of the spin representation of $K'$.}

\ni By the result of lemma~\ref{h.weight}, any highest weight of the spin representation of $K$ has the form
$$\frac{1}{2}\left(\sum_{\beta\in \Phi_2^{+}} \pm \beta'+\sum_{\gamma\in \Phi_3'^{+}} \gamma'\right)\,.$$

\begin{lem} The spin representation of $K'$ has three highest weights both of multiplicity $2$~:
\begin{align*}
\lambda_1&:=7\, \alpha_1'+9\,\alpha_2'+10\,\alpha_3'+5\,\alpha_4'\,,\\
\lambda_2&:=6\, \alpha_1'+9\,\alpha_2'+11\,\alpha_3'+6\,\alpha_4'\,,\\
\lambda_3&:=5\, \alpha_1'+9\,\alpha_2'+12\,\alpha_3'+6\,\alpha_4'\,.\\
\end{align*}
\end{lem}
\proof Note first that
$$\sum_{\gamma\in \Phi_3'^{+}} \gamma'=5\, \alpha_1'+3\,\alpha_2'+6\,\alpha_3'+3\,\alpha_4'\,.$$
Now, denoting by $\beta_1,\ldots,\beta_8$ the roots occurring in that order in the expression of 
$\Phi_2^{+}$ above, it is easily checked that
\begin{align*}
\lambda_1&=
\frac{1}{2}(\beta_1'+\beta_2'+\beta_3'+\beta_4'+
\beta_5'+\beta_6'+\beta_7'+\beta_8')+\frac{1}{2}
\sum_{\gamma\in \Phi_3'^{+}} \gamma'\,,\\
\lambda_2&=\frac{1}{2}(-\beta_1'+\beta_2'+\beta_3'+\beta_4'+
\beta_5'+\beta_6'+\beta_7'+\beta_8')+\frac{1}{2}
\sum_{\gamma\in \Phi_3'^{+}} \gamma'\,,\\
\lambda_3&=\frac{1}{2}(-\beta_1'-\beta_2'+\beta_3'+\beta_4'+
\beta_5'+\beta_6'+\beta_7'+\beta_8')+\frac{1}{2}
\sum_{\gamma\in \Phi_3'^{+}} \gamma'\,,\\
\end{align*}
so $\lambda_1$, $\lambda_2$ and $\lambda_3$ are actually weights of the spin representation
of $K'$.

\ni It is also easily checked that they are also $G$-dominant. Indeed, considering the basis of fundamental weights
$\omega_1'$, $\omega_2'$, $\omega_3'$, $\omega_4'$ associated with $\alpha_1'$, $\alpha_2'$,
$\alpha_3'$, $\alpha_4'$, one gets\footnote{multiplying for instance the transpose of the Cartan matrix with
the vector representing the weight in the basis ($\alpha_1'$, $\alpha_2'$,
$\alpha_3'$, $\alpha_4'$).}
\begin{align*}
\lambda_1&=5\, \omega_1'+\omega_2'+\omega_3'\,,\\
\lambda_2&=3\, \omega_1'+\omega_2'+\omega_3'+\omega_4'\,,\\
\lambda_3&=\omega_1'+\omega_2'+3\, \omega_3'\,.\\
\end{align*}
Note now that, for any $i,j=1,2,3$, $i\not =j$, any relation of the form 
$\lambda_i\prec \lambda_j$, where ``{$\prec$}'' is the standard order on weights, is impossible, 
so $\lambda_i$ can not be an element of the set of weights of an irreducible module with
highest weight $\lambda_j$. This implies that $\lambda_1$, $\lambda_2$ and $\lambda_3$ are actually
highest weights of the spin representation of $K'$.

\ni Since $\dim \gt_0=2$, one knows from \S~\ref{w.spin} that the multiplicity of each such highest weight is at
least $2$. Now, denoting by $n_i$ the dimension of any irreducible module with highest weight $\lambda_i$, 
$i=1,2,3$, one gets with the  help of the LiE program
$$ n_1=180224\,,\quad n_2=524288\quad \text{and}\quad n_3=344064\,.$$
Now 
$$2\times n_1+2\times n_2+2\times n_3=2^{21}=2^{\dim(\gp)/2}=\dim \Sigma\,,$$
hence the result.\qed

\begin{remark} Actually, the determination of the three possible highest weights was obtained by using the same argument as
in \S~\ref{hw.so}. Indeed the $\alpha$'s in $\Phi_1$ and the $\beta$'s in $\Phi_2$
appear to be respectively compact and non-compact roots for the inner symmetric space 
$\mathrm{SO}(8)/(\mathrm{SO}(4)\times\mathrm{SO}(4))$, relatively to the standard torus 
of $\mathrm{SO}(8)$. Using the results of R. Parthasarathy
in \cite{Par}, we obtain, with the help of the GAP3 program, the twelve highest weights 
of the spin representation for the symmetric space $\mathrm{SO}(8)/(\mathrm{SO}(4)\times\mathrm{SO}(4))$,
providing a list of twelve candidates
for the possible highest weights of the spin representation of $K'$. Finally, the dominance condition (for $K'$) reduces
that list to the three weights $\lambda_1$, $\lambda_2$ and $\lambda_3$.
\end{remark}

\subsection{The first eigenvalue of the Dirac operator.}

\ni We determine $G$-weights $\mu$ (non necessarily dominant) such that
$\mu_{|\gt_K}=\lambda_1$, $\lambda_2$ or $\lambda_3$, and $\|\mu+\delta_G\|^2$ is minimal.

\ni Let $\mu=\sum_{i=1}^6 \mu_i\, \theta_i\in i\, \gt^{*}$. First
$$\mu_{|\gt_K}=\mu_2\,\alpha_1'+(\mu_1+\mu_6)\,\alpha_2'+(\mu_3+\mu_5-\mu_2)\,\alpha_3'+
(\mu_4-\mu_2)\,\alpha_4'\,,
$$
so $\mu_{|\gt_K}=\lambda_i$, $i=1,2,3$ if and only if
\begin{align*}
 &\mu_1+\mu_6=9\,,\quad \mu_2=7\,,\quad \mu_3+\mu_5=17\,,\quad \mu_4=12\,,\\
 &\mu_1+\mu_6=9\,,\quad \mu_2=6\,,\quad \mu_3+\mu_5=17\,,\quad \mu_4=12\,,\\
 &\mu_1+\mu_6=9\,,\quad \mu_2=5\,,\quad \mu_3+\mu_5=17\,,\quad \mu_4=11\,.
 \end{align*}

\ni Let $\mu$ be such a weight. We first examine the case where it is $G$-dominant.
\subsubsection{First case: $\mu$ is $G$-dominant.}
The weight $\mu$ is $G$-dominant if and only if 
$$\begin{cases}
 2\,\mu_1-\mu_3\in \bN\,,\\
 \mu_2=7\;\text{and}\; \mu_4=12\,,\;\text{or}\; 
 \mu_2=6\;\text{and}\; \mu_4=12\,,\\
 -\mu_1+2\,\mu_3-12\in \bN\,,\\
 \mu_1-2\,\mu_3+13\in \bN\,,\\
 -2\,\mu_1+\mu_3+1\in \bN\,.\\
\end{cases}$$
Thus setting
$$\begin{cases}k:=2\,\mu_1-\mu_3\,,\\
l:=-\mu_1+2\,\mu_3-12\,,\end{cases}\quad \left(\Longleftrightarrow
\begin{cases}\mu_1=\frac{2k+l+12}{3}\,,\\
\mu_3=\frac{k+2l+24}{3}\,,\end{cases}\right)$$
$$\text{$\mu$ is $G$-dominant}\Longleftrightarrow
\begin{cases} k=0\;\text{or}\; 1\,,\\
l=0\;\text{or}\; 1\,,\\
\mu_2=7\;\text{and}\; \mu_4=12\,,\;\text{or}\; 
 \mu_2=6\;\text{and}\; \mu_4=12\,.
\end{cases}$$

\ni The term involving $\mu_2$ in the expression of 
$\|\mu+\delta_G\|^2$ is 
$$2\big((\mu_2+11)^2-(\mu_2+11)(\mu_4+21)\big)\,,$$
hence the minimum of $\|\mu+\delta_G\|^2$ is obtained only when $\mu_2=6$ (and then $\mu_4=12$).

\ni Viewing then $\|\mu+\delta_G\|^2$ as a polynomial $f(k,l)$ of the variables $k$ and $l$, one gets
$$\frac{\partial f}{\partial k}(k,l)=\frac{2}{3}\,(6\,\mu_1-27)\quad \text{and}\quad 
\frac{\partial f}{\partial l}(k,l)=\frac{2}{3}\,(6\,\mu_3-51)\,.$$
Thus $f(k,l)$ has only one critical point $(\frac{1}{2},\frac{1}{2})$. Now
$$\frac{\partial^2 f}{\partial k^2}(k,l)=\frac{8}{3}\,,\quad 
\frac{\partial^2 f}{\partial k\,\partial l}(k,l)=\frac{4}{3}\quad\text{and}
\quad \frac{\partial^2 f}{\partial l^2}(k,l)=\frac{8}{3}\,,$$
So
$$f(k,l)-f\left(\frac{1}{2},\frac{1}{2}\right)=\frac{2}{3}\,\left((k-\frac{1}{2})^2+
(l-\frac{1}{2})^2+(k+l-1)^2\right)\,.$$
Hence among the dominant $G$-weights $\mu$ such that $\mu_{|\gt_K}=\lambda_i$, $i=1,2,3$, 
the minimum of $\|\mu+\delta_G\|^2$ is obtained
if and only if 
\begin{align*}
 &(k,l)=(0,1)\,, \left(\Leftrightarrow (\mu_1,\mu_3)=\left(\frac{13}{3},\frac{26}{3}\right)\right)\,,\\
 \intertext{or}
 &(k,l)=(1,0)\,, \left(\Leftrightarrow (\mu_1,\mu_3)=\left(\frac{14}{3},\frac{25}{3}\right)\right)\,,
\end{align*}
hence for the two weights 
\begin{align*}
\mu_1&=\frac{13}{3}\,\theta_1+6\, \theta_2+\frac{26}{3}\, \theta_3+12\,\theta_4+\frac{25}{3}\,
\theta_5+\frac{14}{3}\,\theta_6\,,\\
\intertext{or}
\mu_2&=\frac{14}{3}\,\theta_1+6\, \theta_2+\frac{25}{3}\, \theta_3+12\,\theta_4+\frac{26}{3}\,
\theta_5+\frac{13}{3}\,\theta_6\,,
\end{align*}
and in that case $\mu_{|\gt_K}=\lambda_2$. 

\ni It is easily verified that
$$\langle \mu_1+2\delta_G,\mu_1\rangle=\langle \mu_2+2\delta_G,\mu_2\rangle=\frac{340}{3}\,.$$

\subsubsection{Second case: $\mu$ is not $G$-dominant.}
\ni In that case, $k$ and $l$ are arbitrary integers, and, considering the term involving $\mu_2$ in the expression of 
$\|\mu+\delta_G\|^2$, the minimum of $\|\mu+\delta_G\|^2$ is obtained only when $\mu_2=5$ (and then $\mu_4=11$).
Up to a constant term, $\|\mu+\delta_G\|^2$ is then the same polynomial $f(k,l)$ as above, so the minimum 
is obtained for $k=0$ and $l=1$, or $k=1$ and $l=0$, hence we may concude that among the $G$-weights $\mu$ such that
$\mu_{|\gt_K}=\lambda_i$, $i=1,2,3$, 
the minimum of $\|\mu+\delta_G\|^2$ is obtained for the two weights
\begin{align*}
\mu_1'&=\frac{13}{3}\,\theta_1+5\, \theta_2+\frac{26}{3}\, \theta_3+11\,\theta_4+\frac{25}{3}\,
\theta_5+\frac{14}{3}\,\theta_6\,,\\
\intertext{or}
\mu_2'&=\frac{14}{3}\,\theta_1+5\, \theta_2+\frac{25}{3}\, \theta_3+11\,\theta_4+\frac{26}{3}\,
\theta_5+\frac{13}{3}\,\theta_6\,,
\end{align*}
and in that case $\mu_{|\gt_K}=\lambda_3$. 

\ni Now each weight is conjugate under the action of the Weyl group $W_G$ to one and only one 
dominant $G$-weight. Using the function ``{W\_orbit}'' of the LiE program, it may be checked that 
$\mu_1'$ is conjugate to $\mu_1$, and $\mu_2'$ to $\mu_2$. Hence we may conclude that
\begin{lem} Among the $G$-dominant weights $\mu$ verifying the spin condition, 
$$\text{$\|\mu+\delta_G\|^2$ is minimal if and only if $\mu=\mu_1$ or $\mu=\mu_2$.}$$
\end{lem}
Hence the square of the first eigenvalue of the Dirac operator is
$$\frac{1}{24} \langle \mu_1+2\,\delta_G,\mu_1\rangle +\frac{\dim(\gp)}{16}=\frac{85}{18}+\frac{21}{8}=
\frac{529}{72}\,.$$

\section{Appendix}

\subsection{Proof of lemma~\ref{lem.roots}}
\begin{description}
 \item[R1] For any $\theta\in \Phi_3$, one has 
 \begin{align}\label{R11}
 \forall X\in\gt_K\,,\quad&[X,\sigma_{*}(E_{\theta})]=[\sigma_{*}(X),\sigma_{*}(E_{\theta})]=
 \sigma_{*}\left([X,E_{\theta}]\right)=\theta(X)\, \sigma_{*}(E_{\theta})\,,\\
 \intertext{and}\label{R12}
 \forall Y\in\gt_0\,,\quad&[Y,\sigma_{*}(E_{\theta})]=-[\sigma_{*}(Y),\sigma_{*}(E_{\theta})]=
 -\sigma_{*}\left([Y,E_{\theta}]\right)=-\theta(Y)\, \sigma_{*}(E_{\theta})\,.
 \end{align}
Hence 
\begin{align}\label{R13}
\forall X\in\gt_K\,,\; [X,U_{\theta}]=\theta(X)\,U_{\theta}\quad \text{and}\quad [X,V_{\theta}]=\theta(X)\,V_{\theta}\,. 
\end{align}
Now, suppose that for any $X\in \gt_K$, $\theta(X)=0$. Then there exists $Y\in \gt_0$ such that $\theta(Y)\not =0$.
Set $U_\theta=A_\theta+i\, B_\theta$, where $A_\theta$ and $B_\theta\in\gg$.
By \eqref{R13}, for any $X\in \gt_K$, one has $[X,U_{\theta}]=0$, which implies that both $A_\theta$
and $B_\theta$ belong to 
the centralizer of $\gt_K$ in $\gk$ which is equal to $\gt_K$, as $\gt_K$ is maximal. Hence $U_\theta\in 
\gt_{K,\bC}$, so $[Y,U_{\theta}]=0$.
But using \eqref{R12}, $[Y,U_{\theta}]=\theta(Y)\, V_{\theta}$, which can not be $0$ as $\theta(Y)\not =0$ and
$V_{\theta}\not= 0$. Hence, there exists $X\in \gt_K$ such that $\theta(X)\not =0$, and so by \eqref{R13}, 
$\theta'$ is a $K$-root. \qed
\item[R2] Let $\theta\in \Phi_1\cup\Phi_2$. As $\sigma_{*}(E_{\theta})=\pm E_{\theta}$, one has 
for any $Y\in\gt_0$, $[Y,\sigma_{*}(E_{\theta})]=\theta(Y)\,\sigma_{*}(E_{\theta})$. But \eqref{R12} implies 
$[Y,\sigma_{*}(E_{\theta})]=-\theta(Y)\,\sigma_{*}(E_{\theta})$, hence $\theta(Y)=0$.

\ni Conversely, if $\theta_{|\gt_0}=0$, then by \eqref{R11}, for any $X\in \gt$, 
$$[X,\sigma_{*}(E_{\theta})]=\theta(X)\, \sigma_{*}(E_{\theta})\,,$$ so $\sigma_{*}(E_{\theta})\in g_{\theta}$. Hence 
there exists $\lambda\in \bC$ such that $\sigma_{*}(E_{\theta})=\lambda\, E_{\theta}$. As $\sigma_{*}$ is an involution,
one has $\lambda^2=1$, hence $\sigma_{*}(E_{\theta})$ belongs to $\gk_{\bC}$ or $\gp_{\bC}$,
so $\theta \in \Phi_1\cup\Phi_2$.
\item[R3] By \textbf{R1}, $\beta'$ is a $K$-root, and by \textbf{R2}, there exists $Y\in \gt_0$ such that $\beta(Y)\not=0$. Suppose $\beta'=\pm \alpha'$. Then
by \eqref{R13}
$$\forall X\in \gt_K\,, \quad [X,E_{\alpha}]=\alpha(X)\, E_{\alpha}\quad\text{and}\quad
[X,U_{\beta}]=\pm \alpha(X)\, U_{\beta}\,.$$
Hence $U_{\beta}\in \gk_{\pm \alpha'}$, so there exists $\lambda\in \bC$, such that $U_{\beta}=\lambda\, E_{\pm \alpha}$.
Now by \textbf{R2}, on has $[Y,E_{\pm \alpha}]=0$, so $[Y,U_{\beta}]=0$. But by \eqref{R12}, $[Y,U_{\beta}]=\beta(Y)\,
V_{\beta}$, which is nonzero.\qed
\item[R4] Let $\theta\in \Phi_1\cup\Phi_2$. For any $X\in \gt_K$, one has 
$\sigma^{*}(\alpha)(X)=\alpha(\sigma_{*}(X))=\alpha(X)$, and for any $Y\in \gt_0$, by \textbf{R2}, 
$\sigma^{*}(\alpha)(Y)=\alpha(\sigma_{*}(Y))=\alpha(-Y)=0=\alpha(Y)$, hence $\sigma^{*}(\alpha)=\alpha$.

\ni Now let $\theta\in \Phi_3$. By \textbf{R2}, there exists $Y\in \gt_0$ such that $\theta(Y)\not = 0$. Note that 
\begin{align*}
 \forall X\in \gt_K\,,\quad[X,\sigma_{*}(E_{\theta})]&=[\sigma_{*}X,\sigma_{*}(E_{\theta})]=\sigma_{*}\left([X,E_{\theta}]\right)\\
 &=\theta(X)\, \sigma_{*}(E_{\theta})=
 \sigma^{*}(\theta)(X)\, \sigma_{*}(E_{\theta})\,,\\
 \forall X\in \gt_0\,,\quad [X,\sigma_{*}(E_{\theta})]&=-[\sigma_{*}X,\sigma_{*}(E_{\theta})]=-\sigma_{*}\left([X,E_{\theta}]\right)\\
 &=-\theta(X)\, \sigma_{*}(E_{\theta})=
 \sigma^{*}(\theta)(X)\, \sigma_{*}(E_{\theta})\,.
\end{align*}
Hence $\sigma^{*}(\theta)$ is a root. As $\sigma^{*}(\theta)(Y)=-\theta(Y)\not=0$, $\sigma^{*}(\theta)\in \Phi_3$ by \textbf{R2}.
Furthermore, as for any $X\in \gt_K$, $\sigma^{*}(\theta)(X)=\theta(X)$, $\sigma^{*}(\theta)\not =-\theta$, and 
as $\sigma^{*}(\theta)(Y)=-\theta(Y)\not =\theta(Y)$, $\sigma^{*}(\theta)\not =\theta$. 
\item[R5] As $\sigma^{*}$ is an involution, any root $\theta$ may be decomposed into $\theta=\theta_{+}+\theta_{-}$, with
$\sigma^{*}(\theta_{+})=\theta_{+}$ and $\sigma^{*}(\theta_{-})=-\theta_{-}$. Note that $\theta_{+|\gt_0}=0$ and 
$\theta_{-|\gt_K}=0$ since
\begin{align*}
 \forall X\in\gt_0\,,\quad &\theta_{+}(X)=\sigma^{*}(\theta_{+})(X)=\theta_{+}(\sigma_{*}(X))=-\theta_{+}(X)\,,\\
 \intertext{and}
 \forall X\in\gt_K\,,\quad &\theta_{-}(X)=-\sigma^{*}(\theta_{-})(X)=-\theta_{-}(\sigma_{*}(X))=-\theta_{-}(X)\,.
 \end{align*}
This implies in particular that the decomposition is orthogonal, and $\theta_{+}=\theta_{|\gt_K}=\theta'$.

\ni Now let $\alpha$ and $\beta\in \Phi_3$ such that $\beta'=\alpha'$, with $\beta\not=\alpha$ and
$\beta\not=\sigma^{*}(\alpha)$. We first claim that $\alpha+\beta$ is not a root. Otherwise 
$$[E_{\alpha},E_{\beta}]=\underbrace{[U_\alpha,U_\beta]+[V_\alpha,V_\beta]}_{\in \gk_{\bC}}+
\underbrace{[U_\alpha,V_\beta]+[V_\alpha,U_\beta]}_{\in \gp_{\bC}}\,,$$
 is a nonzero element in $\gg_{\alpha+\beta}$.
If $[U_\alpha,U_\beta]+[V_\alpha,V_\beta]=0$, then $[E_{\alpha},E_{\beta}]\in \gp_{\bC}$, so $\alpha+\beta \in \Phi_2$,
hence $(\alpha+\beta)_{|\gt_0}=0$, by \textbf{R2}, so $\beta_{|\gt_0}=-\alpha_{|\gt_0}$, but this is impossible since
$\beta\not=\sigma^{*}(\alpha)$. 

\ni Hence $[U_\alpha,U_\beta]+[V_\alpha,V_\beta]\not =0$. But then $2\alpha'$ is a $K$-root since, as $\beta'=\alpha'$,
$$\forall X\in \gt_K\,, \quad [X,[U_\alpha,U_\beta]+[V_\alpha,V_\beta]]=2\alpha(X)\, 
([U_\alpha,U_\beta]+[V_\alpha,V_\beta])\,.$$
But $\alpha'$ is a $K$-root (\textbf{R1}) and it is well known that $\alpha'$ and $2\alpha'$ can not be both $K$-roots.

\ni Thus $\alpha+\beta$ is not a root, and so $\langle \alpha,\beta\rangle\geq 0$, since otherwise $\alpha+\beta$ should
be a root (see for instance 9.4. in \cite{Hum}).

\ni One the other hand, $\alpha-\beta$ can not be a root by \textbf{R1} since $(\alpha-\beta)_{|\gt_K}=0$. Thus
$\langle \alpha,\beta\rangle\leq 0$,  since otherwise $\alpha-\beta$ should
be a root.

\ni Finally, $\langle \alpha,\beta\rangle=0$. Note that, as $\sigma^{*}(\alpha)_{|\gt_K}=\alpha_{|\gt_K}$,
the same arguments hold for $\sigma^{*}(\alpha)$ instead of $\alpha$. Hence $\langle \sigma^{*}(\alpha),\beta\rangle=0$, 
and so $\langle \alpha+\sigma^{*}(\alpha), \beta\rangle=0$.
Now as $\alpha+\sigma^{*}(\alpha)=2\, \alpha_{+}$, this implies $\langle \alpha_{+},\beta_{+}\rangle=0$.
But $\alpha_{+}=\alpha'=\beta'=\beta_{+}$, hence $\|\alpha'\|^2=0$, which is impossible. \qed
\end{description}

------------------------------------------------

\end{document}